\documentclass[a4paper,reqno]{amsart}
\usepackage{dsfont}
\usepackage{newcent} 
\usepackage[T1]{fontenc}
\usepackage{amsmath,amssymb,amsfonts,amsthm}
\usepackage{mathtools}
\usepackage{newlfont}
\usepackage{fancyhdr,fancyvrb}
\usepackage{graphicx}
\usepackage{nextpage}
\usepackage[dvipsnames,usenames]{color}
\usepackage[colorlinks=true, linkcolor=blue, citecolor=red]{hyperref}
\usepackage{tikz}
\usepackage{tikz-3dplot}
\usepackage{mathrsfs}
\usepackage{stmaryrd}
\usepackage{stackrel}
\usepackage{enumerate} 
\usepackage{enumitem}

\makeatletter
\@namedef{subjclassname@2020}{%
\textup{2020} Mathematics Subject Classification}
\makeatother

\newcommand{\compresslist}{
	\setlength{\itemsep}{1.5pt}
	\setlength{\parskip}{0pt}
	\setlength{\parsep}{0pt}
}

\newcommand{\N}{\mathbb{N}}
\newcommand{\R}{\mathbb{R}}
\renewcommand{\S}{\mathbb{S}}

\newcommand{\Z}[1]{\mathbb{Z}^{#1}}

\newcommand{\cB}{\mathcal B}
\newcommand{\cC}{\mathcal C}

\newcommand{\cH}{\mathcal H}

\newcommand{\cL}{\mathcal L}
\newcommand{\cM}{\mathcal M}
\newcommand{\cO}{\mathcal O}

\newcommand{\cR}{\mathcal R}

\newcommand{\bulk}{\mathrm{bulk}}

\newcommand{\de}{\mathrm d}

\newcommand{\surface}{\mathrm{surf}}

\newcommand{\eps}{\varepsilon}

\newcommand{\dist}[2]{\operatorname{dist}(#1,#2)}

\newcommand{\sgn}{\operatorname{sgn}}

\usepackage{soul}

\renewcommand{\geq}{\geqslant}
\renewcommand{\leq}{\leqslant}

\newcommand{\longrightharpoonup}{\relbar\joinrel\rightharpoonup}
\newcommand{\wto}{\rightharpoonup}

\newcommand{\wsto}{\stackrel{*}{\rightharpoonup}}

\newcommand{\wSD}{\stackrel[SD]{*}{\rightharpoonup}}
\newcommand{\wpSD}{\stackrel[SD_p]{}{\longrightharpoonup}}

\newcommand{\average}{{\mathchoice {\kern1ex\vcenter{\hrule
height.4pt width 8pt depth0pt}
\kern-11pt} {\kern1ex\vcenter{\hrule height.4pt width 4.3pt
depth0pt} \kern-7pt} {} {} }}
\newcommand{\ave}{\average\int}

\newcommand{\res}{\mathop{\hbox{\vrule height 7pt width .5pt depth
0pt\vrule height .5pt width 6pt depth0pt}}\nolimits}

\mathchardef\emptyset="001F

\numberwithin{equation}{section}

\newtheorem{defin}{Definition}[section]
\newtheorem{remark}[defin]{Remark}
\newtheorem{theorem}[defin]{Theorem}
\newtheorem{lemma}[defin]{Lemma}
\newtheorem{proposition}[defin]{Proposition}

\newtheorem{ass}[defin]{Assumptions}

\allowdisplaybreaks

\title[]{Periodic homogenization in the context of structured deformations}%
\author{Micol Amar}
\address[M.~Amar]{Dipartimento di Scienze di Base e Applicate per l'Ingegneria, Sapienza Universit\`{a} di Roma, Via Antonio Scarpa, 16, 00161 Roma, Italy}
\email{micol.amar@uniroma1.it}
\author{Jos\'{e} Matias}
\address[J.~Matias]{Instituto Superior T\'{e}cnico, Universidade de Lisboa, Av.~Rovisco Pais 1, 1049-001 Lisboa, Portugal}
\email{jose.c.matias@tecnico.ulisboa.pt}
\author{Marco Morandotti}
\address[M.~Morandotti]{Dipartimento di Scienze Matematiche ``G.~L.~Lagrange'', Politecnico di Torino, Corso Duca degli Abruzzi, 24, 10129 Torino, Italy}
\email{marco.morandotti@polito.it}
\author{Elvira Zappale}
\address[E.~Zappale]{Dipartimento di Scienze di Base e Applicate per l'Ingegneria, Sapienza Universit\`{a} di Roma, Via Antonio Scarpa, 16, 00161 Roma, Italy}
\email{elvira.zappale@uniroma1.it}

\date{\today}


\makeindex

\begin{document}

\begin{abstract}
An energy for first-order structured deformations in the context of periodic homogenization is obtained. This energy, defined in principle by relaxation of an initial energy of integral type featuring contributions of bulk and interfacial terms, is proved to possess an integral representation in terms of relaxed bulk and interfacial energy densities. These energy densities, in turn, are obtained via asymptotic cell formulae defined by suitably averaging, over larger and larger cubes, the bulk and surface contributions of the initial energy. The integral representation theorem, the main result of this paper, is obtained by mixing blow-up techniques, typical in the context of structured deformations, with the averaging process proper of the theory of homogenization.
\end{abstract}

\maketitle

\tableofcontents


\section{Introduction}\label{sec_intro}

The mathematical modeling of materials has interested scientists for many centuries. In the last decades, the accuracy of these models has considerably increased as an effect of more sophisticated measure instruments, on the experimental side, and of the availability of sound mathematical abstract frameworks, on the theoretical side.
Classical theories of continuum mechanics provide a good description of many phenomena such as elasticity, plasticity, and fracture, and are susceptible of incorporating fine structures at the microscopic level. The mathematical process through which the effects of the microstructure emerge at the macroscopic level is called homogenization: this procedure provides an effective macroscopic description as the result of averaging out the heterogeneities.

Structured deformations \cite{DPO1993} provide a mathematical framework to capture the effects at the macroscopic level of geometrical changes at submacroscopic levels. The availability of this framework, especially in its variational formulation \cite{CF1997}, leads naturally to the enrichment of the energies and force systems that underlie variational and field-theoretic descriptions of important physical phenomena without having to commit at the outset to any of the existing prototypical mechanical theories, such as elasticity or plasticity.
A (first-order) structured deformation is a pair $(g,G)\in SBV(\Omega;\R^d)\times L^1(\Omega;\R^{d\times N})\eqqcolon SD(\Omega)$, where $g\colon\Omega\to\R^d$ is the macroscopic deformation and $G\colon\Omega\to\R^{d\times N}$ is the microscopic deformation tensor. 
As opposed to classical theories of mechanics, in which~$g$ and its gradient~$\nabla g$ alone characterize the deformations of the body~$\Omega$, the additional geometrical field~$G$ captures the contributions at the macroscopic level of the smooth submacroscopic changes. The difference $\nabla g-G$ captures the contributions at the macroscopic level of slips and separations occurring at the submacroscopic level (which are commonly referred to as \emph{disarrangements} \cite{DPO1993}). 
Heuristically, the disarrangement tensor $M\coloneqq \nabla g-G$ is an indication of how non classical a structured deformation is: should~$M=0$ and if~$g$ is a Sobolev field, then the field~$G$ is simply the classical deformation gradient; on the contrary, if $M\neq0$, there is a macroscopic bulk effect of submacroscopic slips and separations, which are phenomena involving interfaces. 
This fact will be made precise in the Approximation Theorem~\ref{app_thm}.

In order to assign an energy to a structured deformation $(g,G)\in SD(\Omega)$, the proposal has been made in \cite{CF1997} to take the energetically most economical way to reach $(g,G)$ by means of $SBV$ fields $u_n$: according to the Approximation Theorem \cite[Theorem~2.12]{CF1997}, we say that a sequence $\{u_n\}\subset SBV(\Omega;\R^d)$ converges to $(g,G)$ if
\begin{equation}\label{101}
u_n\to g\quad\text{in $L^1(\Omega;\R^d)$}\qquad\text{and}\qquad \nabla u_n\wsto G\quad\text{in $\cM(\Omega;\R^{d\times N})$,}
\end{equation}
where $\cM(\Omega;\R^{d\times N})$ is the set of bounded matrix-valued Radon measures on $\Omega$, and we denote this convergence by $u_n\wSD (g,G)$.
We let the initial energy of a deformation $u\in SBV(\Omega;\R^d)$ be
\begin{equation}\label{103}
E(u)\coloneqq \int_\Omega W(\nabla u(x))\,\de x+\int_{\Omega\cap S_u} \psi([u](x),\nu_u(x))\,\de\cH^{N-1}(x),
\end{equation}
which is determined by the bulk and surface energy densities $W\colon\R^{d\times N}\to[0,+\infty)$ and $\psi\colon\R^d\times\S^{N-1}\to[0,+\infty)$.  In formula~\eqref{103}, $\de x$ and $\de\cH^{N-1}(x)$ denote the $N$-dimensional Lebesgue and $(N-1)$-dimensional Hausdorff measures, respectively; $[u](x)$ and $v_u(x)$ denote the jump of $u$ and the normal to the jump set for each $x\in S_u$, the jump set.

In mathematical terms, the process just described to assign the energy to a structured deformation $(g,G)\in SD(\Omega)$ 
reads
\begin{equation}\label{102}
I(g,G)\coloneqq \inf\Big\{\liminf_{n\to\infty} E(u_n): \{u_n\}\subset SBV(\Omega;\R^d),\,u_n\wSD(g,G)\Big\}.
\end{equation}
In the language of calculus of variations, the operation described in \eqref{102} is called \emph{relaxation} and the main results in \cite{CF1997} was to prove that the functional $I$ admits an integral representation, that is, there exist functions $H\colon\R^{d\times N}\times\R^{d\times N}\to[0,+\infty)$ and $h\colon\R^d\times\S^{N-1}\to[0,+\infty)$ such that
\begin{equation}\label{104}
I(g,G)=\int_\Omega H(\nabla g(x),G(x))\,\de x+\int_{\Omega\cap S_g} h([g](x),\nu_g(x))\,\de\cH^{N-1}(x).
\end{equation}

In this work, we focus on heterogeneous, hyperelastic, defective materials featuring a fine periodic microstructure. Our scope is to provide an asymptotic analysis of the energies associated with these materials, as the fineness of their microstructure vanishes, in the variational context of structured deformations \cite{CF1997}.

The initial energy functionals that we consider involve a bulk contribution and a surface contribution, each of which is described by an energy density which depends explicitly on the spatial variable in a periodic fashion, namely the energy associated with a deformation $u\in SBV(\Omega;\R^d)$ has the expression
\begin{equation}\label{100}
E_\eps(u)\coloneqq \int_\Omega W\Big(\frac{x}{\eps},\nabla u(x)\Big)\de x+\int_{\Omega\cap S_u} \psi\Big(\frac{x}{\eps},[u](x),\nu_u(x)\Big)\de\cH^{N-1}(x),
\end{equation}
where $W\colon \R^N\times\R^{d\times N}\to[0,+\infty)$ and $\psi\colon\R^N\times\R^d\times\S^{N-1}\to[0,+\infty)$ are $Q$-periodic in the first variable ($Q$~being the unit cube in $\R^N$), and $\eps>0$ is the length scale of the microscopic heterogeneities (see Assumptions~\ref{ass} for the precise assumptions on $W$ and $\psi$).
We will perform a relaxation analogous to that in \eqref{102} for the $\eps$-dependent initial energy \eqref{100}.
In particular, we aim at assigning an energy to structured deformation $(g,G)$ where the geometric field $G$ is $p$-integrable for some $p>1$. As proved in \cite{CF1997}, this has the effect that the submacroscopic slips and separations diffuse in the bulk and contribute to determining the relaxed bulk energy density, whereas the relaxed surface energy density is only determined by optimizing the initial surface energy density.
The mechanical interpretation of this fact, which is also the motivation for our choice, is that the relaxed surface energy is not influenced by the elastic energy in the limit.

We define the class of admissible sequences for the relaxation by
\begin{equation}\label{105-1}
\cR_p(g,G;\Omega)\coloneqq\Big\{\{u_n\}\in SBV(\Omega;\R^d): u_n\wSD(g,G) ,\, \sup_{n\in\N}\lVert \nabla u_n\rVert_{L^p(\Omega;\R^{d\times N})}<+\infty
\Big\}
\end{equation}
and for every sequence $\eps_n\to0$, we define 
\begin{equation}\label{105}
I_{\hom}(g,G)\!\coloneqq \inf\!\Big\{\!\liminf_{n\to\infty} E_{\eps_n}(u_n): \{u_n\}\in\cR_p(g,G;\Omega)\Big\}.
\end{equation}
The main result of this work is the following theorem, which provides a representation result analogous to that in \eqref{104}.
Let $SD_p(\Omega)\coloneqq SBV(\Omega;\R^d)\times L^p(\Omega;\R^{d\times N})$.
\begin{theorem}\label{thm_main}
Let $p>1$ and let us assume that Assumptions~\ref{ass} hold; let $u\in SBV(\Omega;\R^d)$ and let $E_\eps(u)$ be the energy defined by \eqref{100}. 
Then, for every $(g,G)\in SD(\Omega)$, the homogenized functional $I_{\hom}(g,G)$ defined in \eqref{105} admits the integral representation 
\begin{equation}\label{106}
I_{\hom}(g,G)=\begin{cases}
\displaystyle \int_\Omega H_{\hom}(\nabla g(x),G(x))\,\de x+ \int_{\Omega\cap S_g} h_{\hom}([g](x),\nu_g(x))\,\de\cH^{N-1}(x) \\[3mm]
\hfill \text{if $(g,G)\in SD_p(\Omega)$,} 
\\[2mm]
+\infty \hfill \text{otherwise.}
\end{cases}
\end{equation}
The relaxed energy densities $H_{\hom}\colon\R^{d\times N}\times\R^{d\times N}\to[0,+\infty)$ and $h_{\hom}\colon\R^d\times\S^{N-1}\to[0,+\infty)$ are given by the formulae 
\begin{equation}\label{bulkd}
	\begin{split}
\!\!\!\!		H_{\hom} (A, B) \coloneqq \inf_{k \in \N{}} \frac{1}{k^N} \inf \bigg\{&\, \int_{kQ} \!\! W(x, A + \nabla u(x))\, \de x  \\
		& + \int_{kQ\cap S_u} \!\!\!\!  \psi (x, [u](x),\nu_u(x))\, \de \cH^{N-1}(x): u\in\cC_p^{\bulk}(A,B;kQ)\bigg\}
\end{split}\end{equation}
for every $A,B\in\R^{d\times N}$,
and 
\begin{equation}\label{108}
\!\! h_{\hom}(\lambda,\nu)\coloneqq 
\inf_{k\in\N} 
\frac{1}{k^{N-1}}\inf\bigg\{ 
\int_{(kQ_\nu)\cap S_u}  \!\!\!\! \!\!\!\! \psi(x,[u](x),\nu_u(x))\,\de\cH^{N-1}(x): 
 u\in \cC^{\surface}(\lambda,\nu;k
Q_\nu)\bigg\}  
\end{equation}
for every $(\lambda,\nu)\in\R^d\times\S^{N-1}$, where $Q_\nu$ is any rotated unit cube so that two faces are perpendicular to $\nu$. \end{theorem}

The independence of $h_{\hom}(\lambda,\nu)$ from the specific choice of the cube $Q_\nu$ can be deduced from Proposition~\ref{gammahomcont}.
In \eqref{bulkd} and \eqref{108}, we have defined, for $A,B\in\R^{d\times N}$, $(\lambda,\nu)\in\R^d\times\S^{N-1}$, and $R,R_\nu\subset\R^N$ cubes,
\begin{eqnarray}
\cC_p^{\bulk}(A,B;R) &\!\!\!\!\coloneqq&\!\!\!\! \bigg\{u \in SBV_{\#}(R;\R^d): |\nabla u| \in L^p(R), \ave_{R} \nabla u(x)\,\de x = B-A \bigg\}, \label{Cpbulk} \\
\!\!\!\!\!\!\!\! \!\!\!\!\!\!\!\!\!\!\!\!\!\!\!\! \cC^{\surface}(\lambda,\nu;R_\nu) &\!\!\!\!\coloneqq&\!\!\!\! \big\{ u\in SBV(R_\nu;\R^d): u|_{\partial R_\nu}(x)=s_{\lambda,\nu}(x),  \nabla u(x)=0\; \text{a.e.~in $R_\nu$}\big\}, \label{Csurface} 
\end{eqnarray}
where
\begin{equation}\label{special}
s_{\lambda,\nu}(x)=\frac12\lambda(\sgn(x\cdot\nu)+1)
\end{equation}
is 
the elementary jump of amplitude~$\lambda$ across the hyperplane perpendicular to~$\nu$. 
In formula \eqref{Cpbulk}, we denote by $SBV_{\#}(R;\R^d)$ the set of $\R^d$-valued $SBV$ functions with equal traces on opposite faces of the cube $R$.

We notice that~\eqref{bulkd} and~\eqref{108} are asymptotic cell formulae, as it is expected in the context of homogenization when no convexity assumptions are made on the initial energy densities (see, \emph{e.g.}, \cite{BD1998}).
In the special case of functions $W$ and $\psi$ which are convex in the gradient and jump variable, respectively, we are able to show that \eqref{bulkd} reduces to a cell problem in the unit cell (see Proposition~\ref{Hhomprop} below); whether the same result holds for $h_{\hom}$ is still unknown.

Since, in a structured deformation $(g,G)\in SD(\Omega)$, the field $G$ is generally different from~$\nabla g$, it is clear that the convergence \eqref{101} is obtained at the expenses of the discontinuity sets of $u_n$ diffusing in the bulk, namely $\cH^{N-1}(S_{u_n})\to+\infty$ as $n\to\infty$ (so that the hypotheses of Ambrosio's compactness theorem in $SBV$ \cite{A1989} are in general not satisfied). 
This is reflected in the form of the relaxed bulk energy density in \eqref{bulkd}, where we point out that both the initial bulk and surface energy densities contribute to the definition of $H_{\hom}$, both undergoing the bulk rescaling.
On the contrary, the coercivity assumption (see Assumption~\ref{ass}-\ref{H8} below) yields an $L^p$ constraint 
on the gradients of the approximating sequences which avoids the appearance of any bulk contributions in the relaxed surface energy density $h_{\hom}$ in \eqref{108}.

The proof of \eqref{106}, to which the whole Section~\ref{sec_main_proof} is devoted, is obtained by computing the $\Gamma$-limit in \eqref{105} by combining blow-up techniques \emph{\`{a} la} Fonseca-M\"{u}ller \cite{FM1992,FM1993} with rescaling techniques typically used in homogenization problems.
This will be especially visible in the construction of the recovery sequences for proving that the densities in \eqref{106} are indeed given by~\eqref{bulkd} and~\eqref{108}.
To deduce the upper bound for the homogenized surface energy density $h_{\hom}$ we also make use of comparison results in a $\Gamma$-convergence setting (see \cite{BDV1996,DalMaso}).

We will collect some preliminary results in Section~\ref{sec_prel}, where we also prove the Approximation Theorem~\ref{app_thm} which guarantees the non-emptiness of the class $\cR_p$ introduced in \eqref{105-1}. 
Section~\ref{sec_main} contains the precise formulation of the standing assumptions on the initial energy densities $W$ and $\psi$ and a collection of results on the homogenized energy densities $H_{\hom}$ and $h_{\hom}$ which can be deduced from the definitions \eqref{bulkd} and \eqref{108}.
For the reader's convenience, we present in Appendix~\ref{appa} some technical measure-theoretical results which are by now standard.

\section{Preliminaries}\label{sec_prel}
\subsection{Notation}
We will use the following notations
\begin{itemize}
\item $\N$ denotes the set of natural numbers without the zero element;
\item $\Omega \subset \R^{N}$ is a bounded connected open set with Lipschitz boundary;
\item $\S^{N-1}$ denotes the unit sphere in $\R^N$;
\item  for any $r>0$, $B_r$ denotes the open ball of $\R^{N}$ centred at the origin of radius~$r$; for any $x\in\R^{N}$, $B_r(x) \coloneqq x+ B_r$ denotes the open ball centred at $x$ of radius~$r$; $Q\coloneqq (-\tfrac12,\tfrac12)^N$ denotes the open unit cube of $\R^{N}$ centred at the origin; for any $\nu\in\S^{N-1}$, $Q_\nu$ denotes any open unit cube in $\R^{N}$ with two faces orthogonal to~$\nu$; for any $x\in\R^{N}$ and $\delta>0$, $Q(x,\delta)\coloneqq x+\delta Q$ denotes the open cube in $\R^{N}$ centred at $x$ with side $\delta$;
\item ${\mathcal A}(\Omega)$ is the family of all open subsets of $\Omega $; 
\item $\cL^{N}$ and $\cH^{N-1}$ denote the  $N$-dimensional Lebesgue measure and the $\left(  N-1\right)$-dimensional Hausdorff measure in $\R^N$, respectively; the symbol $\de x$ will also be used to denote integration with respect to $\cL^{N}$; 
\item $\cM(\Omega;\R^{d\times N})$ is the sets of finite matrix-valued Radon measures on $\Omega$; $\cM ^+(\Omega)$ is the set of non-negative finite Radon measures on $\Omega$;
given $\mu\in\cM(\Omega;\R^{d\times N})$, 
the measure $|\mu|\in\cM^+(\Omega)$ 
denotes the total variation of $\mu$;
\item $SBV(\Omega;\R^d)$ is the set of vector-valued \emph{special functions of bounded variations} defined on $\Omega$. 
Given $u\in SBV(\Omega;\R^d)$, its distributional gradient $Du$ admits the decomposition $Du=D^au+D^su=\nabla u\cL^N+[u]\otimes\nu_u\cH^{N-1}\res S_u$, where $S_u$ is the jump set of~$u$, $[u]$ denotes the jump of~$u$ on $S_u$, and $\nu_u$ is the unit normal vector to $S_u$; finally, $\otimes$ denotes the dyadic product; 
for $Q\subset\R^N$ a cube, we denote by $SBV_{\#}(Q;\R^d)$ the set of $\R^d$-valued $SBV$ functions with equal traces on opposite faces of $Q$;
\item $L^p(\Omega;\R^{d\times N})$ is the set of matrix-valued $p$-integrable functions; for $p>1$ we denote by $p'$ its H\"{o}lder conjugate.
\item for $p\geq1$, $SD_p(\Omega)\coloneqq SBV(\Omega;\R^d)\times L^p(\Omega;\R^{d\times N})$ is the space of structured deformations $(g,G)$ (notice that $SD_1(\Omega)$ is the space $SD(\Omega)$ introduced in \cite{CF1997}); 
\item  $C$ represents a generic positive constant that may change from line to line.
\item For every $x \in \mathbb R^N$, the symbol $\lfloor x\rfloor \in \mathbb Z^N$ denotes the integer part of the vector~$x$, namely that vector whose components are the integer parts of each component of~$x$. We denote by $ \langle x \rangle$ the fractional part of $x$, i.e., $\langle x \rangle\coloneqq x-\lfloor x\rfloor\in[0,1)^N$.
\end{itemize}

\subsection{Function spaces}
The following proposition serves as a definition of Lebesgue points for $L^p$ functions (see \cite[Theorem~1.33]{EG2015} for a more general statement).
\begin{proposition}[Lebesgue points]\label{propLP}
Let $p\geq1$ and let $u\in L^p(\Omega)$. Then for $\cL^N$-a.e.~$x_0\in\Omega$, the following equality holds
\begin{equation}\label{LP}
\lim_{r\to0^+} \ave_{Q(x_0;r)} |u(x)-u(x_0)|^p\,\de x=0.
\end{equation}

\end{proposition}

The following theorem collects some facts about $BV$ functions. Its proof can be found, e.g., in \cite[Sections~3.6 and~3.7]{AFP2000}, \cite[Section~6.1]{EG2015}, and \cite[Theorem~4.5.9]{Federer1969}.
\begin{theorem}\label{thm2.3BBBF}
Let $u\in BV(\Omega;\mathbb{R}^{d})$. Then
\begin{enumerate}
\item[(i)] (approximate differentiability) for $\mathcal{L}^{N}$-a.e.~$x_0\in\Omega$%
\begin{equation*}
\lim_{r\to0^{+}}\frac{1}{r}\bigg\{\ave_{Q(x_0;r)}\vert u(x)  -u(x_0)  -\nabla u(x_0)\cdot(x-x_0)\vert ^{\frac{N}{N-1}}\,\de x\bigg\}  ^{\frac{N-1}{N}}=0; 
\end{equation*}
\item[(ii)] (jump points) for every $x_0\in S_u$\,, there exist $u^{+}(x_0), u^{-}(x_0)\in\mathbb{R}^{d}$ and $\nu(x_0)  \in \S^{N-1}$ normal to $S_u$ at $x_0$ such that
$$\lim_{r\to0^{+}}\frac{1}{r^{N}}\int_{Q_{\nu(x_0)}^{\pm}(x_0;r)}\big\vert u(x)  -u^{\pm}(x_0)\big\vert\, \de x=0, 
$$
where $Q_{\nu(x_0)}^{\pm}(x_0;r) \coloneqq \{x\in Q_{\nu(x_0)}(x_0;r) : (x-x_0) \cdot \nu(x_0) \gtrless0\}$;
\item[(iii)] (Lebesgue points) for $\cH^{N-1}$-a.e.~$x_0\in\Omega\setminus S_u$\,, \eqref{LP} holds true.
\end{enumerate}
\end{theorem}

Observe that $(i)$ above entails
\begin{equation}\label{appdiffasCF}
\lim_{r\to0^{+}}\frac{1}{r^{N+1}}\int_{Q(x_0;r)}\vert u(x)  -u(x_0)  -\nabla u(x_0)\cdot(x-x_0)\vert \,\de x =0\end{equation}

\subsection{The approximation theorem in $SD_p(\Omega)$}
In this section we prove the approximation theorem for structured deformations in $SD_p(\Omega)$. This result will be useful for the proof of our homogenization Theorem~\ref{thm_main} and rests on the following two statements.
\begin{theorem}[{\cite[Theorem~3]{A1991}}]\label{Al}
Let $f \in L^1(\Omega; \R^{d\times N})$. 
Then there exist $v \in SBV(\Omega; \R^d)$ and a Borel function $\beta\colon\Omega\to\R^{d\times N}$ such that
\begin{equation}\label{817}
Dv = f \cL^N + \beta \cH^{N-1}\res S_v, \quad
\int_{\Omega\cap S_v} |\beta(x)| \, \de \cH^{N-1}(x) \leq C_N \lVert f\rVert_{L^1(\Omega; \R^{d\times N})\,,}
\end{equation}
where $C_N>0$ is a constant depending only on $N$. 
\end{theorem}
\begin{theorem}[{\cite[Lemma~2.9]{CF1997}}]\label{ctap}
Let $v \in BV(\Omega; \R^d)$. Then there exist piecewise constant functions $\bar v_n\in SBV(\Omega;\R^d)$  such that $\bar v_n \to v$ in $L^1(\Omega; \R^d)$ and
\begin{equation}\label{818}
|Dv|(\Omega) = \lim_{n\to \infty}| D\bar v_n|(\Omega) = \lim_{n\to \infty} \int_{\Omega\cap S_{\bar v_n}} |[\bar v_n](x)|\, \de\cH^{N-1}(x).
\end{equation}
\end{theorem}
One of the main results in the theory developed by Del Piero and Owen was the Approximation Theorem, stating that any structured deformation can be approximated, in the $L^\infty$ sense, by a sequence of simple deformations (see \cite{DPO1993} for the details, in particular Theorem~5.8).
For structured deformations $(g,G)\in SD(\Omega)$, the corresponding result is obtained in \cite[Theorem~2.12]{CF1997}.
Here we prove a version in $SD_p(\Omega)$, which is the natural framework for the integral representation of the functional $I_{\hom}$ defined in \eqref{105}.
\begin{theorem}[Approximation Theorem]\label{app_thm}
For every $(g,G)\in SD_p(\Omega)$ there exists a sequence $u_n\in SBV(\Omega;\R^d)$ such that $u_n\wpSD (g,G)$, namely
\begin{equation}\label{110}
u_n\to g\quad\text{in $L^1(\Omega;\R^d)$}\qquad\text{and}\qquad \nabla u_n\wto G\quad\text{in $L^p(\Omega;\R^{d\times N})$.}\end{equation}
Moreover, there exists $C>0$ such that, for all $n\in\N{}$,
\begin{equation}\label{appEST}
|D u_n|(\Omega)\leq C \big(\|g\|_{BV(\Omega;\R^{d})}+\|G\|_{L^p(\Omega;\R^{d\times N})}\big). 
\end{equation}
In particular, 
this implies that, up to a subsequence, 
\begin{equation}\label{centerline}
D^s u_n\wsto (\nabla g-G)\cL^N+D^s g\qquad \text{in $\cM(\Omega;\R^{d\times N})$.}
\end{equation}
\end{theorem}
\begin{proof}
Let $(g,G)\in SD_p(\Omega)$ and, by Theorem~\ref{Al} with $f\coloneqq \nabla g-G$, let $v\in SBV(\Omega;\R^{d})$ be such that $\nabla v=\nabla g-G$.
Furthermore, let $\bar v_n\in SBV(\Omega;\R^{d})$ be a sequence of piecewise constant functions approximating $v$, as per Lemma~\ref{ctap}.
Then, the sequence of functions 
\begin{equation*}
u_n\coloneqq g+\bar v_n-v
\end{equation*}
is easily seen to approximate $(g,G)$ in the sense of \eqref{110}. 
In fact, $ u_n \to g$ in $L^1(\Omega;\R^d)$ and $\nabla u_n(x) = G(x)$ for $\cL^N$-a.e.~$x\in\Omega$.
Estimate \eqref{appEST} follows from the inequality in \eqref{817} and from \eqref{818}; 
finally, \eqref{110} and \eqref{appEST} imply \eqref{centerline}.
\end{proof}


\section{Standing assumptions and properties of the homogenized densities}\label{sec_main}
In this section we present the hypotheses on the initial energy densities $W$ and $\psi$ and we prove some properties of the homogenized densities $H_{\hom}$ and $h_{\hom}$ defined in \eqref{bulkd} and \eqref{108}, respectively.

\begin{ass}\label{ass}
Let $p>1$ and let $W\colon \R^N \times \R^{d\times N} \to[0,+\infty)$ and $\psi\colon \R^N\times \R^d \times \S^{N-1}\to[0,+\infty)$ be continuous functions such that 
\begin{enumerate}
\item \label{H1} for every $\xi\in\R^{d\times N}$ and for every $(\lambda,\nu)\in\R^d\times\S^{N-1}$, the functions $x\mapsto W(x,\xi)$ and $x\mapsto \psi(x, \lambda, \nu)$ are $Q$-periodic;
\item \label{H2} there exists $C_W > 0$ such that, for every $x\in\R^N$ and for every $\xi_1,\xi_2\in\R^{d\times N}$, 
$$|W(x,\xi_1) - W(x,\xi_2)| \leq C_W|\xi_1-\xi_2|( 1 + |\xi_1|^{p-1} + |\xi_2|^{p-1});$$
\item \label{H2c} there exists a function $\omega_W\colon [0,+\infty)\to[0,+\infty)$ such that $\omega_W(s)\to0$ as $s\to0^+$ such that for every $x_1,x_2\in\R^{N}$ and $\xi\in\R^{d\times N}$
$$|W(x_1,\xi)-W(x_2,\xi)|\leq\omega_W(|x_1-x_2|)(1+|\xi|^p);$$
\item \label{H8} there exist $C'_W>0$, and $c'_W>0$ such that $W(x,\xi)\geq C'_W |\xi|^p- c'_W$ for every $\xi \in \mathbb R^{d\times N}$ and a.e.~$x\in \Omega$.
\item \label{H3} there exists $c_\psi, C_\psi > 0$ such that, for every $(x,\lambda,\nu)\in\R^N\times\R^d\times\S^{N-1}$,
$$c_\psi|\lambda| \leq \psi(x, \lambda, \nu) \leq C_\psi|\lambda|;$$
\item \label{H3c}  there exists a function $\omega_\psi\colon [0,+\infty)\to[0,+\infty)$ such that $\omega_\psi(s)\to0$ as $s\to0^+$ such that for every $x_1,x_2\in\R^{N}$ and $(\lambda,\nu)\in\R^{d}\times\S^{N-1}$
$$|\psi(x_1,\lambda,\nu)-\psi(x_2,\lambda,\nu)|\leq\omega_\psi(|x_1-x_2|)|\lambda|;$$ 
\item \label{H4} for every $(x,\nu)\in\R^N\times\S^{N-1}$, the function $\lambda\mapsto\psi(x,\lambda,\nu)$ is \emph{positively homogeneous of degree one}, \emph{i.e.}, for every $\lambda\in\R^d$ and $t>0$,
$$\psi(x, t\lambda, \nu) = t\psi(x, \lambda, \nu);$$
\item \label{H5} for every $(x,\nu)\in\R^N\times\S^{N-1}$, the function $\lambda\mapsto\psi(x,\lambda,\nu)$ is \emph{subadditive},\emph{i.e.}, for every $\lambda_1,\lambda_2\in\R^d$,
$$ \psi(x, \lambda_1 + \lambda_2, \nu) \leq \psi(x, \lambda_1, \nu) + \psi(x, \lambda_2, \nu);$$
\item \label{H6} for every $x\in\R^N$, the function $(\lambda,\nu)\mapsto\psi(x,\lambda,\nu)$ is \emph{symmetric}, \emph{i.e.}, for every $(\lambda,\nu)\in\R^d\times\S^{N-1}$,
$$ \psi(x, \lambda, \nu) = \psi(x, - \lambda, - \nu);$$
\end{enumerate}
\end{ass}
\begin{remark}\label{rem1}
We make the following observations.
\begin{enumerate}
\item The $p$-Lipschitz continuity in \ref{H2} jointly with \ref{H1} and \ref{H2c} imply  that $W$ has $p$-growth from above in the second variable, namely that there exists $C_W>0$ such that for every $(x,\xi)\in\R^N\times\R^{d\times N}$
\begin{equation}\label{pgrabove}
W(x,\xi)\leq C_{W}(1+|\xi|^p).
\end{equation}
\color{black}
On the contrary, $p$-growth from above jointly with the quasiconvexity of the bulk energy density in the gradient variable (which is the natural assumption in equilibrium problems in elasticity) returns the $p$-Lipschitz continuity.
\item Condition \ref{H3} does not allow for a control on the $\cH^{N-1}$-measure of the jump set, which, in the spirit of the Approximation Theorem~\ref{app_thm}, is crucial in the context of structured deformations.
\item Conditions \ref{H3} and \ref{H5} imply Lipschitz continuity of the function $\lambda\mapsto\psi(x,\lambda,\nu)$, \emph{i.e.}, for every $(x,\nu)\in\R^N\times\S^{N-1}$ and for every $\lambda_1,\lambda_2\in\R^d$,
\begin{equation}\label{Lippsi}
|\psi(x,\lambda_1, \nu)- \psi(x,\lambda_2, \nu)| \leq C_\psi|\lambda_1- \lambda_2|,
\end{equation}
\item Conditions \ref{H4}, \ref{H5}, and \ref{H6} are natural ones for fractured materials; in particular, condition \ref{H6} allows one to identify $\psi(x,\lambda,\nu)$ with $\widetilde\psi(x,\lambda\otimes\nu)$, for a suitable function $\widetilde\psi\colon\R^N\times\R^{d\times N}\to[0,+\infty)$.
\end{enumerate}
\end{remark}



We now present a translation invariance property of $H_{\hom}$ and $h_{\hom}$. 
\begin{proposition}[translation invariance]\label{118a}
For $A,B\in\R^{d\times N}$, let $H_{\hom}(A,B)$ be defined by \eqref{bulkd}. Then for every $\tau\in Q$, we have $H_{\hom} (A, B) = H^\tau_{\hom}(A,B)$, where 
\begin{equation}\label{118b}
\begin{split}
H^\tau_{\hom} (A, B) \coloneqq \inf_{k \in \N} \frac{1}{k^N}&\, \inf \bigg\{ \int_{kQ} \!\! W(x+\tau, A + \nabla u(x))\, \de x \\
&\, + \int_{kQ\cap S_u} \!\!\!\! \psi (x+\tau, [u](x),\nu_u(x))\, \de \cH^{N-1}(x): 
u\in\cC_p^{\bulk}(A,B;kQ)
\bigg\},
\end{split}
\end{equation}
where $\cC_p^{\bulk}(A,B;kQ)$ is defined in \eqref{Cpbulk}.\\
For $(\lambda,\nu)\in\R^d\times\S^{N-1}$, let $h_{\hom}(\lambda,\nu)$ be defined by \eqref{108}. Then for every $\tau\in Q$, we have $h_{\hom} (\lambda,\nu) = h^\tau_{\hom}(\lambda,\nu)$, where 
\begin{equation}\label{108tau}
\begin{split}
h_{\hom}^\tau(\lambda,\nu)\coloneqq 
\inf_{k\in\N} 
\frac{1}{k^{N-1}}\inf\bigg\{ & 
\int_{(kQ_\nu)\cap S_u}  \psi(x+\tau,[u](x),\nu_u(x))\,\de\cH^{N-1}(x): \\
&\, u\in \cC^{\surface}(\lambda,\nu;k
Q_\nu)\bigg\},  
\end{split}
\end{equation}
where $\cC^{\surface}(\lambda,\nu;kQ_\nu)$ is defined in \eqref{Csurface}.
\end{proposition}
\begin{proof}
The proof of both \eqref{118b} and \eqref{108tau} is a straightforward adaptation of the proof of \cite[Proposition 2.15]{MMS}.
\end{proof}


The next proposition shows that under if the initial bulk and surface energy densities $W$ and $\psi$ are convex in the gradient and jump variable, respectively, then the homogenized bulk energy density is obtained via a cell formula in in the unit cube and is not asymptotic anymore. 
\begin{proposition}\label{Hhomprop}
Let $W$ and $\psi$ satisfy Assumptions~\ref{ass}, let us assume that the functions $\xi\mapsto W(x,\xi)$ and $\lambda\mapsto\psi(x,\lambda,\nu)$ are convex for every $x\in\R^d$ and every $\nu\in\S^{N-1}$, and let
\begin{equation*}
	\begin{split}
H_{\hom}^{\mathrm{cell}} (A, B) \coloneqq  \inf \bigg\{&\, \int_{Q} \!\! W(x, A + \nabla u(x))\, \de x  \\
		& + \int_{Q\cap S_u}  \psi (x, [u](x),\nu_u(x))\, \de \cH^{N-1}(x): u\in\cC_p^{\bulk}(A,B;Q)\bigg\}.
\end{split}
\end{equation*}
Then $H_{\hom}(A,B)=H_{\hom}^{\mathrm{cell}}(A,B)$ for every $A,B\in\R^{d\times N}$.
\end{proposition}
\begin{proof}
Let $A,B\in\R^{d\times N}$ be given and let us denote by $m_k(A,B)$ the inner infimization problem in the definition of $H_{\hom}(A,B)$, so that \eqref{bulkd} reads $H_{\hom}(A,B)=\inf_{k\in\N} m_k(A,B)$. With this position, we also have $H_{\hom}^{\mathrm{cell}}(A,B)=m_1(A,B)$.

We obtain the desired result if we prove that $m_k(A,B)=m_1(A,B)$. 
To this aim, let $u\in\cC_p^{\bulk}(A,B;Q)$ be an admissible function for $m_1(A,B)$. By extending $u$ by $Q$-periodicity on $kQ$, we obtain a function in $\cC_p^{\bulk}(A,B;kQ)$ which is a competitor for $m_k(A,B)$, whence $m_k(A,B)\leq m_1(A,B)$. 
To show the reverse inequality, we consider $u\in\cC_p^{\bulk}(A,B;kQ)$ a competitor for $m_k(A,B)$ and we use the standard method of averaging its translates to produce a competitor $v\in\cC_p^{\bulk}(A,B;Q)$ for $m_1(A,B)$, see \cite[proof of Theorem~14.7]{BD1998}, and using Jensen's inequality.
By letting $J\coloneqq\{0,1,\ldots,k-1\}^N$, it is easy to see that the function $v\colon \R^{N}\to\R^{d}$ defined by
$$Q\ni x\mapsto v(x)\coloneqq \frac1{k^N}\sum_{j\in J}u(x+j)$$
and extended by periodicity is $Q$-periodic and satisfies 
\[\begin{split}
\int_Q \nabla v(x)\,\de x=&\, \frac1{k^N}\sum_{j\in J} \int_Q \nabla u(x+j)\,\de x=\frac1{k^N} \sum_{j\in J} \int_{Q-j} \nabla u(y)\,\de y \\
=&\, \frac1{k^N} \int_{\bigcup_{j\in J} (Q-j)} \nabla u(y)\,\de y = \ave_{kQ} \nabla u(y)\,\de y=B-A,
\end{split}\]
so that $v\in\cC_p^{\bulk}(A,B;Q)$ and therefore $m_1(A,B)\leq m_k(A,B)$, yielding the sought-after equality $m_k(A,B)=m_1(A,B)$ and the independence of the size of the cube. The thesis follows.
\end{proof}

The next proposition contains further properties of $h_{\hom}$.

\begin{proposition}\label{gammahomcont}
Let $\psi$ satisfy Assumptions~\ref{ass} and let $\widehat h_{\hom}\colon \R^d\times\S^{N-1}\to[0,+\infty)$ be the function defined by
\begin{equation}\label{108hat}
\begin{split}
\widehat h_{\hom}(\lambda,\nu)\coloneqq \limsup_{T\to+\infty} \frac1{T^{N-1}} \inf\bigg\{ & \int_{(TQ_\nu)\cap S_u} \psi(x,[u](x),\nu_u(x))\,\de\cH^{N-1}(x): \\
&\, u\in \cC^{\surface}(\lambda,\nu;TQ_\nu)\bigg\}.  
\end{split}
\end{equation}
Then the following properties hold true:
\begin{enumerate}
\item \label{ha}
the function $\widehat h_{\hom}$ 
is a limit which is independent of the choice of the cube $Q_\nu$, once $\nu\in\S^{N-1}$ is fixed;
\item \label{hb} the function $\widehat h_{\hom}$ is continuous on $\R^d\times\S^{N-1}$ and for every $(\lambda, \nu) \in \R^d\times \S^{N-1}$ 
\begin{equation}\label{gammahomgrowth}
c_{\psi}|\lambda|\leq  \widehat h_{\hom}(\lambda,\nu) \leq 
C_{\psi}|\lambda|,
\end{equation}
where $c_\psi$ and $C_\psi$ are the constants in Assumptions~\ref{ass}-\ref{H3};
\item \label{hc} for every $(\lambda,\nu)\in\R^d\times \S^{N-1}$, we have 
$h_{\hom}(\lambda,\nu)=\widehat h_{\hom}(\lambda,\nu)$,
where $h_{\hom}$ is the function defined in \eqref{108}.
\end{enumerate}
\end{proposition}
\begin{proof}
The proof of items \ref{ha} and \ref{hb} is essentially the same as that of \cite[Proposition 2.2]{BDV1996},
upon observing that our density $\psi$ satisfies \eqref{Lippsi}, which is a stronger continuity assumption than condition \cite[(iii) page~304]{BDV1996}.
Conclusion \ref{ha} is obtained \emph{verbatim} as in \cite[proof of Proposition~2.2, Steps 1--4]{BDV1996}; we sketch here a proof of conclusion \ref{hb} for the reader's convenience.

The continuity of $\widehat h_{\hom}$ can be obtained by arguing in the following way:
\begin{itemize}\compresslist
\item[(a)] one shows that the function $\widehat h_{\hom}(\lambda,\cdot)$ is continuous on $\S^{N-1}$, uniformly with respect to $\lambda$, when $\lambda$ 
varies on bounded sets; 
\item[(b)] one shows that for every $\nu\in \S^{N-1}$, the function $\widehat h_{\hom}(\cdot,\nu)$ is continuous on $\R^d$;
\item[(c)] one shows that $\widehat h_{\hom}$ is continuous in the pair $(\lambda,\nu)$.
\end{itemize}
The proof of point (a) above relies on the fact that for every fixed $\lambda\in\R^d$, formula \eqref{108hat} does not depend on the cube $Q_\nu$ once the direction $\nu$ is prescribed, by \ref{ha}.
For the proof of point (b),
we can argue as in \cite[proof of Proposition~2.2, Step 6]{BDV1996} (here we exploit Assumptions~\ref{ass}-\ref{H3} and the Lipschitz continuity of $\psi$, see \eqref{Lippsi}).
Point (c) can be obtained by arguing as in \cite[proving (ii) from (i) in Theorem~2.8]{DGpont}.

To conclude the proof of \ref{hb}, we need to prove \eqref{gammahomgrowth}.
The estimate from above can be easily obtained from the very definition of $\widehat h_{\hom}$ in \eqref{108hat}, by using Assumptions~\ref{ass}-\ref{H3}.
Concerning the estimate from below, it is sufficient to observe that the functional $SBV(\Omega;\R^d)\ni u\mapsto \int_{S_u}|[u](x)|\,\de\cH^{N-1}(x)$ is lower semicontinuous with respect to the convergence $u_n \wpSD(g, 0)$, with $g$ a pure jump function, as it follows from the lower semicontinuity of the total variation with respect to the weak-* convergence and, again, from Assumptions~\ref{ass}-\ref{H3}.

To prove \ref{hc}, we take inspiration from the proof of \cite[Lemma~2.1]{CDDA2006}: we show that $\widehat h_{\hom}$ is an infimum over the integers. Together with \ref{ha}, we will conclude that $\widehat h_{\hom}=h_{\hom}$, as desired. 
Let $g_T\colon\R^d\times\S^{N-1}\to[0,+\infty)$ be defined by
\begin{equation}\label{g_T}
\!\!\! g_T(\lambda,\nu)\coloneqq  \frac1{T^{N-1}} \inf\bigg\{ \int_{(TQ_\nu)\cap S_u} \!\!\! \psi(x,[u](x),\nu_u(x))\,\de\cH^{N-1}(x):\, u\in \cC^{\surface}(\lambda,\nu;TQ_\nu)\bigg\},
\end{equation}
so that we can write \eqref{108hat} as $\widehat h_{\hom}(\lambda,\nu)=\limsup_{T\to+\infty} g_T(\lambda,\nu)$.

We start by proving a monotonicity property of $g_T$ over multiples of integer values of $T$, namely we prove that, for every $(\lambda,\nu)\in\R^d\times\S^{N-1}$, 
\begin{equation}\label{108m}
g_{hk}(\lambda,\nu)\leq g_k(\lambda,\nu)\qquad\text{for every $h,k\in\N
$.}
\end{equation}
To this aim, let $u\in\cC^{\surface}(\lambda,\nu;kQ_\nu)$ be a competitor for $g_k(\lambda,\nu)$ and consider the function $\bar u\colon hkQ_\nu\to\R^d$ defined by
$$\bar u(x)\coloneqq\begin{cases}
0 & \text{if $x\cdot\nu<-k/2$,} \\
u(k\langle x/k\rangle) & \text{if $|x\cdot\nu|<k/2$,} \\
\lambda & \text{if $x\cdot\nu>k/2$,}
\end{cases}$$
obtained by replicating $u$ by periodicity in the $(N-1)$-dimensional strip perpendicular to $\nu$ and extending it to $0$ and $\lambda$ appropriately.
It is immediate to see that $\bar u\in\cC^{\surface}(\lambda,\nu;hkQ_\nu)$, so that \eqref{108m} follows.

We now consider two integers $0<m<n$ and a function $u\in\cC^{\surface}(\lambda,\nu,mQ_\nu)$; we define $\tilde u\colon nQ_\nu\to\R^d$ by
$$\tilde u(x)\coloneqq\begin{cases}
u(x) & \text{if $x\in mQ_\nu$,} \\
s_{\lambda,\nu}(x) & \text{if $x\in nQ_\nu\setminus mQ_\nu$}
\end{cases}$$
and notice that $\tilde u\in\cC^{\surface}(\lambda,\nu;nQ_\nu)$.
Then, invoking Assumptions~\ref{ass}-\ref{H3},
\[\begin{split}
&\int_{(nQ_\nu)\cap S_{\tilde u}} \psi(x,[\tilde u](x),\nu_{\tilde u}(x))\,\de\cH^{N-1}(x) \\
=& \int_{(mQ_\nu)\cap S_u} \psi(x,[u](x),\nu_u(x))\,\de\cH^{N-1}(x) + \int_{(nQ_\nu\setminus mQ_\nu)\cap S_{\tilde u}}  \psi(x,[\tilde u](x),\nu_{\tilde u}(x))\,\de\cH^{N-1}(x) \\
\leq& \int_{(mQ_\nu)\cap S_u} \psi(x,[u](x),\nu_u(x))\,\de\cH^{N-1}(x) +C_\psi|\lambda|\big(n^{N-1}-m^{N-1}\big),
\end{split}\]
so that, by infimizing first over $\tilde u\in\cC^{\surface}(\lambda,\nu;nQ_\nu)$ and then over $u\in\cC^{\surface}(\lambda,\nu;mQ_\nu)$, we obtain
$$g_n(\lambda,\nu)\leq g_m(\lambda,\nu) +\frac{C_\psi|\lambda|\big(n^{N-1}-m^{N-1}\big)}{n^{N-1}}.$$
Using $[n/m]m$ in place of $m$ and \eqref{108m}, we get
$$g_n(\lambda,\nu)\leq g_m(\lambda,\nu) +\frac{\displaystyle C_\psi|\lambda|\Big(n^{N-1}-\Big[\frac{n}{m}\Big]^{N-1}m^{N-1}\Big)}{n^{N-1}}.$$
By \ref{ha}, $\widehat h_{\hom}(\lambda,\nu)=\lim_{T\to+\infty} g_T(\lambda,\nu)$, so that, by taking the limit as $n\to\infty$ in the inequality above, we can write
$$\widehat h_{\hom}(\lambda,\nu)=\lim_{T\to+\infty} g_T(\lambda,\nu)=\lim_{n\to\infty} g_n(\lambda,\nu)\leq g_m(\lambda,\nu),\quad\text{for every $m\in\N
$};$$
this yields, by taking the infimum over the integers,
$$\inf_{n\in\N} g_n(\lambda,\nu)\leq \lim_{n\to\infty} g_n(\lambda,\nu)\leq \inf_{m\in\N} g_m(\lambda,\nu),$$
the first inequality being obvious. Recalling the definition \eqref{108} of $h_{\hom}(\lambda,\nu)$, this gives the equality $\widehat h_{\hom}=h_{\hom}$ of \ref{hc} and concludes the proof.
\end{proof} 


\section{Proof of Theorem~\ref{thm_main}}\label{sec_main_proof}
This section is entirely devoted to the proof of Theorem~\ref{thm_main}.
The proof is achieved by obtaining upper and lower bounds for the Radon--Nikod\'ym derivatives of the functional $I_{\hom}$ defined in \eqref{105} with respect to the Lebesgue measure $\cL^N$ and to the Hausdorff measure $\cH^{N-1}$ in terms of the homogenized bulk and surface energy densities $H_{\hom}$ and $h_{\hom}$ defined in \eqref{bulkd} and \eqref{108}, respectively.

We start by observing that any $(g,G)\in SD(\Omega)$ for which there exists a sequence $\{u_n\}\in\cR_p(g,G)$ is indeed an element of $SD_p(\Omega)$. In particular, the functional $(g,G)\mapsto I_{\hom}(g,G)$ is finite if and only if $(g,G)\in SD_p(\Omega)$. Therefore, we have $I_{\hom}\colon SD(\Omega)\to[0,+\infty]$ defined by
\begin{equation}\label{Ihomnew}
I_{\hom}(g,G)\coloneqq 
\begin{cases}
\tilde I_{\hom}(g,G) & \text{if $(g,G)\in SD_p(\Omega)$,} \\
+\infty & \text{if $(g,G)\in SD(\Omega)\setminus SD_p(\Omega)$,}
\end{cases}
\end{equation}
where
\begin{equation}\label{Ihomtilde}
\tilde I_{\hom}(g,G)\coloneqq \inf\Big\{\liminf_{n\to\infty} E_{\eps_n}(u_n): \{u_n\}\in \widetilde\cR_p(g,G;\Omega)\Big\},
\end{equation}
with
\begin{equation}\label{105-2}
\widetilde\cR_p(g,G;\Omega)\coloneqq\Big\{\{u_n\}\in SBV(\Omega;\R^d): u_n\wpSD(g,G)\Big\}.
\end{equation}
We notice that by the Approximation Theorem~\ref{app_thm}, the class $\widetilde\cR_p$ defined above is not empty.
Since if $(g,G)\in SD(\Omega)\setminus SD_p(\Omega)$ there is nothing to prove, we will drop the tildes in the proof.

\subsection{The bulk energy density}\label{sect_bed}
We tackle here the bulk energy density $H_{\hom}$.
In the next two subsections, we assume that  $x_0\in\Omega$ is a point of approximate differentiability for $g$ and a Lebesgue point for $G$, namely, Theorem~\ref{thm2.3BBBF}(i) and (iii) hold for $g$ and \eqref{LP} holds for $G$ (notice that $\cL^N$-a.e.~$x_0\in\Omega$ satisfies these properties).

\subsubsection*{The bulk energy density: lower bound}
Let $\{u_n\}\in \cR_p(g,G;\Omega)$, 
and let $\mu_n\in\cM^+(\Omega)$ be the Radon measure defined by
$$\mu_n \coloneqq W\Big( \frac{x}{\eps_n}, \nabla u_n(x)\Big)\cL^N + \psi \Big( \frac{x}{\eps_n}, [u_n](x), \nu_u(x)\Big)\cH^{N-1}\res S_{u_n}.$$
Without loss of generality, we can assume that $\sup_{n\in\N} \mu_n(\Omega) < + \infty$, so that there exists $\mu\in\cM^+(\Omega)$ such that (up to a not relabelled subsequence) $\mu_n \wsto \mu$.

We will prove that 
\begin{equation}\label{BEDLB}
\frac{\de\mu}{\de \cL^N} (x_0) \geq H_{\hom} (\nabla g(x_0),G(x_0)).
\end{equation}

Let $\{r_k\}$ be a vanishing sequence or radii such that
$\mu(\partial Q(x_0; r_k)) = 0$; then we have
\begin{equation}\label{cube}
\begin{split}
\frac{\de\mu}{\de \cL^N}(x_0) =&  \lim_{k\to \infty} \frac{\mu(Q(x_0; r_k))}{|Q(x_0; r_k)|}= \lim_{k\to\infty}\frac1{r_k^N}\lim_{n\to\infty} \bigg(\int_{Q(x_0; r_k)} W\Big( \frac{x}{\eps_n}, \nabla u_n(x)\Big)\, \de x \\
&  \phantom{\lim_{k\to \infty} \frac{\mu(Q(x_0; r_k))}{|Q(x_0; r_k)|}=}+ \int _{Q(x_0; r_k)\cap S_{u_n}} \psi\Big( \frac{x}{\eps_n}, [u_n](x), \nu_{u_n}(x)\Big)\, \de \cH^{N-1}(x)\bigg) \\
=& \lim_{k\to\infty}\lim_{n\to\infty} \frac1{r_k^N} \bigg(r_k^N \int_Q W\Big( \frac{x_0 + r_ky}{\eps_n}, \nabla_x u_n(x_0 + r_ky)\Big)\, \de y \\
&+ r_k^{N-1}\int_{Q\cap \frac{S_{u_n}- x_0}{r_k}} \psi \Big( \frac{x_0 + r_ky}{\eps_n}, [u_n](x_0 + r_ky), \nu_{u_n}(x_0 + r_ky)\Big)\, \de \cH^{N-1}(y)\bigg),
\end{split}
\end{equation}
where we have changed variables in the last equality.
Upon defining, for every $y\in Q$,
\begin{equation}\label {linear}
u_0(y) \coloneqq \nabla g(x_0)y
\end{equation}
and \begin{equation}\label{unk}
u_{n,k}(y) \coloneqq \frac {u_n(x_0 + r_ky) - g(x_0)}{r_k} - u_0(y),
\end{equation}
we have that 
\begin{equation}\label{863}
\lim_{k\to\infty}\lim_{n\to\infty} \int_Q \big|u_{n,k}(y)\big|\,\de y= 0,
\end{equation}
\begin{equation}\label{924}
\nabla u_{n,k}(y) = \nabla_x u_n (x_0 + r_ky) - \nabla g(x_0),
\end{equation}
and
\begin{equation}\label{weakLp}
	\lim_{k\to\infty}\lim_{n\to\infty} \int_Q (\nabla u_{n,k}(y)- G(x_0)+ \nabla g(x_0))\varphi (y)\de y=0
	\end{equation}
(for any $\varphi \in L^{p'}(Q;\mathbb R^{d \times N}$), where we used the fact that $u_n \wpSD (g,G)$ and Theorems \ref{thm2.3BBBF} and \ref{propLP}. Thus, \eqref{cube} becomes
\begin{equation}\label{cube2}
\begin{split}
\!\!\!\! \frac{\de\mu}{\de \cL^N}(x_0)  =& \lim_{k\to\infty} \lim_{n\to\infty} \frac{1}{r_k^N}\bigg(r_k^N\int_Q W\Big( \frac{x_0 + r_k y}{\eps_n}, \nabla u_{n,k}(y) + \nabla g(x_0)\Big)\, \de y \\
& \phantom{\lim_{k\to\infty} \lim_{n\to\infty} \frac{1}{r_k^N}\bigg(} + r_k^{N-1} \int_{Q \cap S_{u_{n,k}}} \!\!\!\! \!\!  \psi\Big( \frac{x_0 + r_ky}{\eps_n}, r_k [u_{n,k}](y), \nu_{u_{n,k}}(y)\Big)\, \de \cH^{N-1}(y)\bigg)\\
=& \lim_{k\to\infty} \lim_{n\to\infty} \bigg(\int_Q W\Big( \frac{x_0 + r_k y}{\eps_n}, \nabla u_{n,k}(y) + \nabla g(x_0)\Big)\, \de y \\
& \phantom{\lim_{k\to\infty} \lim_{n\to\infty} \bigg(} +\int_{Q \cap S_{u_{n,k}}} \!\!\!\! \!\!  \psi\Big( \frac{x_0 + r_ky}{\eps_n}, [u_{n,k}](y), \nu_{u_{n,k}}(y)\Big)\, \de \cH^{N-1}(y)\bigg),
\end{split}
\end{equation}
where we have used the positive $1$-homogeneity of $\psi$ (see \ref{H4}). 
By writing 
$$ \frac{x_0 + r_ky}{\eps_n} = \frac{r_k y}{\eps_n} + \Big\lfloor \frac{x_0}{\eps_n}\Big\rfloor + \Big \langle \frac{x_0}{\eps_n}\Big \rangle,$$
and using the 1-periodicity of $W$ and $\psi$ (see \ref{H1}), \eqref{cube2} becomes
\begin{equation}\label{cube3}
\begin{split}
\frac{\de\mu}{\de \cL^N}(x_0) = \lim_{k\to\infty}  \lim_{n\to\infty} \bigg(& \int_Q W\Big( \frac{r_ky}{\eps_n} + \Big \langle \frac{x_0}{\eps_n}\Big \rangle, \nabla u_{n,k}(y) + \nabla g(x_0)\Big)\, \de y\\
& + \int_{Q \cap S_{u_{n,k}}} \psi \Big( \frac{r_k y}{\eps_n} + \Big \langle \frac{x_0}{\eps_n}\Big \rangle, [u_{n,k}](y), \nu_{u_{n,k}}(y)\Big)\, \de \cH^{N-1}(y)\bigg).
\end{split}
\end{equation}
We now choose $n(k)$ so that, setting $s_k \coloneqq r_k/\eps_{n(k)}$, we have that $\lim_{k\to\infty} s_k = +\infty$. 
By defining $v_k (y) \coloneqq u_{n(k),k}(y)$ for every $y\in Q$,  \eqref{cube3} becomes
\begin{equation}\label{cube4}
\begin{split}
\frac{\de\mu}{\de \cL^N}(x_0) = \lim_{k\to\infty}  \bigg(& \int_Q W( s_ky + \gamma_k, \nabla v_k(y) + \nabla g(x_0)) \, \de y \\
&+ \int_{Q \cap S_{v_k}} \psi( s_ky + \gamma_k, [v_k](y), \nu_{v_k}(y))\,\de \cH^{N-1}(y)\bigg),
\end{split}
\end{equation}
where $\gamma_k\coloneqq \langle x_0/\eps_{n(k)}\rangle$.
By \eqref{924} we have that $\nabla v_k(y)=\nabla_x u_{n(k)} (x_0 + r_k y) - \nabla g(x_0)$; from \eqref{863} and \eqref{weakLp} 
the sequence $\{v_k\}$ satisfies
\begin{equation}\label{vuk}
v_k \to 0\;\; \text{in $L^1(Q;\R^d)$}\quad\text{and}\quad \nabla v_k  \wto  G(x_0) - \nabla g(x_0)\;\;\text{in $L^p(Q;\R^{d\times N})$}\quad\text{as $k \to \infty$.}
\end{equation}

It is now possible\footnote{
This is achieved, following the strategy in \cite[Proposition~3.1 Step~2]{CF1997}, by constructing a double-indexed sequence that gradually makes a transition from $v_k$ to its limit. The transition takes place across a suitably located frame of vanishing thickness $1/m$ (independent of $k$) and is obtained via convex combination (see also the construction in \cite[Lemma~2.21]{CF1997} where suitable truncations of the approximating sequences are considered; \cite[Lemma~2.20]{CF1997} (see Lemma~\ref{lemma2.20CF}) states that it is possible to work on bounded sequences). A sequence $\tilde w_k$ is then obtained by a diagonalization argument. Finally, the condition on the average in \eqref{961} is enforced by a further modification of the sequence $\tilde w_k$ into $w_k$ by modification with a linear function on cubes that invade $Q$. The difference between our problem and that in \cite{CF1997} is the explicit dependence on the spatial variable that we have; nonetheless, our assumptions \ref{H2} and \ref{H3} allow us to estimate the vanishing terms independently of the spatial variable.} to replace the sequence $\{v_k\}$ with a sequence $\{w_k\}\subset SBV(Q;\R^d)$ still satisfying the convergences in \eqref{vuk}, such that 
\begin{equation}\label{961}
w_k|_{\partial Q}=0\quad\text{and}\quad \int_Q \nabla w_k(y)\,\de y=\frac{(\lfloor s_k\rfloor+1)^N}{s_k^N}(G(x_0)-\nabla g(x_0)) \quad \text{for every $k\in\N$,}
\end{equation} 
and such that
\begin{align*}
\lim_{k\to\infty}  \bigg(\! \int_Q W( s_ky + \gamma_k, \nabla v_k(y) + \nabla g(x_0)) \, \de y +\!\! \int_{Q \cap S_{v_k}} \!\!\!\! \!\! \psi( s_ky + \gamma_k, [v_k](y), \nu_{v_k}(y))\,\de \cH^{N-1}(y)\bigg)\!\geq\\ 
\limsup_{k\to\infty}  \bigg(\! \int_Q W( s_ky + \gamma_k, \nabla w_k(y) + \nabla g(x_0)) \, \de y+\!\! \int_{Q \cap S_{w_k}} \!\!\!\! \!\! \psi( s_ky + \gamma_k, [w_k](y), \nu_{w_k}(y))\,\de \cH^{N-1}(y)\bigg),
\end{align*}
so that \eqref{cube4} becomes
\begin{equation}\label{cube5}
\begin{split}
\frac{\de\mu}{\de \cL^N}(x_0) \geq \liminf_{k\to\infty} \bigg(& \int_Q W( s_ky + \gamma_k, \nabla w_k(y) + \nabla g(x_0))\, \de y \\
&\, + \int_{Q\cap S_{w_k}} \psi (s_k y+ \gamma_k, [w_k](y), \nu_{w_k}(y))\, \de \cH^{N-1}(y)\bigg).
\end{split}
\end{equation}
By changing variables, setting $z\coloneqq s_ky$ and $U_k(z)\coloneqq s_kw_k(x/s_k)$, so that
\begin{equation}\label{694}
U_k|_{\partial(s_kQ)}=0,\quad \nabla_z U_k(z) = \nabla_y w_k\bigg(\frac{z}{s_k}\bigg), \quad\text{and}\quad \frac{1}{s_k}[U_k](z) = [w_k](z),
\end{equation}
we obtain
%
\begin{equation}\label{skcube2}
\begin{split}
\frac{\de\mu}{\de \cL^N}(x_0)\geq \liminf_{k\to\infty} \frac{1}{s_k^N} \bigg(& \int_{s_kQ} W( z + \gamma_k, \nabla U_k(z)+ \nabla g(x_0))\, \de z \\
& + \int_{(s_kQ) \cap S_{U_k}} \psi ( z + \gamma_k, [U_k](z), \nu_{U_k}(z))\, \de \cH^{N-1}(z)\bigg),
\end{split}
\end{equation}
where we have used the positive $1$-homogeneity of $\psi$ once again (see \ref{H4}).

In order to comply with the definition of $H_{\hom}(\nabla g(x_0), G(x_0))$ (see \eqref{bulkd}), we need to integrate over integer multiples of $Q$. To this aim, we extend $U_k$ to the cube $(\lfloor s_k\rfloor+1)Q$ by setting
\begin{equation}\label{hatUk}
\hat U_k(z)\coloneqq\begin{cases}
U_k(z) & \text{if $z\in s_kQ$,} \\[1mm]
0 & \text{if $z\in (\lfloor s_k\rfloor+1)Q\setminus(s_kQ)$.}
\end{cases}
\end{equation}
Notice that, by the first condition in \eqref{694} no further jumps are created, so that  $[\hat U_k](z)=[U_k](z)$ for every $z\in S_{\hat U_k}=S_{U_k}$. Moreover, if follows from \eqref{hatUk}, the definition of $U_k$ and the second condition in \eqref{961} that 
$$
\hat U_k|_{\partial (\lfloor s_k\rfloor+1)Q}=0\quad\text{and}\quad \ave_{(\lfloor s_k\rfloor+1)Q} \nabla\hat U_k(z)\,\de z=G(x_0)-\nabla g(x_0) \quad\text{for every $k\in\N$,}
$$
so that $\{\hat U_k\}\in \cC_p^{\bulk}\big(\nabla g(x_0),G(x_0);(\lfloor s_k\rfloor+1)Q\big)$ (see \eqref{Cpbulk}).
Then, using \eqref{pgrabove}
and the linear growth of $\psi$ (see \ref{H3}), we can continue with \eqref{skcube2} and obtain
\begin{equation}\label{skcube3}
\begin{split}
\frac{\de\mu}{\de \cL^N}(x_0) \geq&\, \liminf_{k\to\infty} \frac{1}{s_k^N} \bigg( \int_{(\lfloor s_k\rfloor + 1)Q} W( z + \gamma_k, \nabla \hat U_k(z)+ \nabla g(x_0))\, \de z \\
&\, \phantom{ \liminf_{k\to\infty} \frac{1}{s_k^N} \bigg(}+ \int_{(\lfloor s_k\rfloor  + 1)Q \cap S_{\hat U_k}} \psi ( z + \gamma_k, [\hat U_k](z), \nu_{\hat U_k}(z))\, \de \cH^{N-1}(z)\bigg)\\
&- \limsup_{k\to +\infty} \frac{1}{s_k^N} \int_{(\lfloor s_k\rfloor + 1)Q\setminus s_k Q} W( z + \gamma_k, \nabla g(x_0))\, \de z  \\ \geq &\liminf_{k\to\infty} \frac{1}{s_k^N} \bigg( \int_{(\lfloor s_k\rfloor + 1)Q} W( z + \gamma_k, \nabla \hat U_k(z)+ \nabla g(x_0))\, \de z \\
&\, \phantom{ \liminf_{k\to\infty} \frac{1}{s_k^N} \bigg(}+ \int_{(\lfloor s_k\rfloor  + 1)Q \cap S_{\hat U_k}} \psi ( z + \gamma_k, [\hat U_k](z), \nu_{\hat U_k}(z))\, \de \cH^{N-1}(z)\bigg)\\
& 
-\limsup_{k \to +\infty}\frac{C_W}{s_k^N}(1+|\nabla g(x_0)|^p) \cL^N(\big(\lfloor s_k\rfloor + 1)Q\setminus s_k Q\big) \\
\geq&\, \liminf_{k\to\infty} \frac{1}{(\lfloor s_k\rfloor+1)^N} \bigg( \int_{(\lfloor s_k\rfloor + 1)Q} W( z + \gamma_k, \nabla \hat U_k(z)+ \nabla g(x_0))\, \de z \\
&\, \phantom{ \liminf_{k\to\infty} \frac{1}{s_k^N} \bigg(}+ \int_{(\lfloor s_k\rfloor  + 1)Q \cap S_{\hat U_k}} \psi ( z + \gamma_k, [\hat U_k](z), \nu_{\hat U_k}(z))\, \de \cH^{N-1}(z)\bigg)\\
\geq&\, \liminf_{k\to\infty} H_{\hom}^{\gamma_k}(\nabla g(x_0), G(x_0)) = H_{\hom}(\nabla g(x_0), G(x_0)),
\end{split}
\end{equation}
where we have used Proposition~\ref{118a} for the last equality.
\qed

\subsubsection*{The bulk energy density: upper bound}
Here we prove that 
\begin{equation}\label{ubb00}
\frac{\de I_{\hom}(g,G)}{\de \cL^N}(x_0)\leq H_{\hom} (\nabla g(x_0),G(x_0)).
\end{equation}
Let $k\in\N\setminus\{0\}$  and $u\in \cC_p^{\bulk}(\nabla g(x_0),G(x_0);kQ)$ (see \eqref{Cpbulk}). 
Let us consider a sequence of radii $r_j\to0$ as $j\to\infty$, and let $h_j\in SBV(Q_{rj}(x_0);\R^d)$ be the function provided by Theorem~\ref{Al} such that 
\begin{equation}\label{ubb03}
\nabla h_j(x)=\nabla g(x_0)-\nabla g(x)+G(x)-G(x_0);
\end{equation}
finally, let $\{h_{j,n}\}_n$ be a piecewise constant approximation of $h_j$ in $L^1(Q_{r_j}(x_0);\R^d)$ provided by Theorem~\ref{ctap}.
We notice that, thanks to Proposition~\ref{propLP} and Theorem~\ref{thm2.3BBBF}, 
\begin{equation}\label{ubb09}
\lim_{j\to\infty} \frac{\alpha_j}{r_j^N}=0,
\end{equation}
where $\alpha_j \coloneqq C\big( \big\lVert |G-G(x_0)|^p \big\rVert_{L^1(Q_{r_j}(x_0))} + \lVert \nabla g - \nabla g(x_0)\rVert_{L^1(Q_{r_j}(x_0);\R^d)}\big)$.
For every $j,n\in\N{}$, we define the function $u_{j,n}\in SBV(Q_{r_j}(x_0);\R^d)$ by
\begin{equation}\label{ubb04}
u_{j,n}(x)\coloneqq g(x)+\frac{r_j}{m_n k}u\Big(\frac{m_n k}{r_j}(x-x_0)\Big)+h_j(x)-h_{j,n}(x),
\end{equation}
where $\{m_n\}_n$ is a diverging sequence of integer numbers to be defined later. 
By defining $kQ\ni y\coloneqq k(x-x_0)/r_j$, and by applying the Riemann--Lebesgue Lemma to the sequence of functions $kQ\ni y\mapsto u^{(n)}(y)\coloneqq u(m_ny)$, we obtain that $u^{(n)}$ converges weakly in $L^p(kQ;\R^d)$ to $\ave_{kQ} u(y)\,\de y$, so that 
\begin{equation}\label{ubb05}
\lim_{n\to\infty} u_{j,n}=g\qquad\text{in $L^1(Q_{r_j}(x_0);\R^d)$} \quad\text{for every $j\in\N{}$;}
\end{equation}
moreover, recalling \eqref{ubb03},  we have
\begin{equation}\label{ubb06}
\begin{split}
\nabla u_{j,n}(x) =&\, \nabla g(x)+\nabla u\Big(\frac{m_n k}{r_y}(x-x_0)\Big)+\nabla h_j(x) \\
=&\, \nabla u\Big(\frac{m_n k}{r_j}(x-x_0)\Big)+\nabla g(x_0)+G(x)-G(x_0),
\end{split}
\end{equation}
so that, by applying the Riemann--Lebesgue Lemma to the sequence $kQ\ni y\mapsto \nabla u^{(n)}(y)\coloneqq \nabla u(m_ny)$, we obtain that $\nabla u^{(n)}$ converges weakly in $L^p(kQ;\R^{d\times N})$ to $\ave _{kQ} \nabla u(y)\,\de y$, yielding
\begin{equation}\label{ubb07}
\lim_{n\to\infty} \nabla u_{j,n}=\ave_{kQ} \nabla u(y)\,\de y+\nabla g(x_0)+G-G(x_0)=G \quad \text{weakly in $L^p(kQ;\R^{d\times N})$.}
\end{equation}
The convergences in \eqref{ubb05} and \eqref{ubb07} show that the sequence $\{u_{j,n}\}_n$
is admissible for the definition of $I_{\hom}(g,G;Q_{r_j}(x_0))$, for every $j\in\N{}$.  

Recalling that the localization $\cO(\Omega)\ni A\mapsto \tilde I_{\hom}(g,G;A)$ is the trace of a Radon measure on the open subsets of $\Omega$ (see Proposition~\ref{asprop2.22CF}) and that for every $(g,G)\in SD_p(\Omega)$ we have that $I_{\hom}(g,G)=\tilde I_{\hom}(g,G)$, we can 
estimate
\begin{equation*}
\begin{split}
&\,\frac{\de I_{\hom}(g,G)}{\de \cL^N}(x_0) \leq \limsup_{j\to\infty}\frac1{r_j^N}\liminf_{n\to\infty} \bigg(\int_{Q_{r_j}(x_0)} \!\! W\Big(\frac{x}{\eps_n},\nabla u_{j,n}(x)\Big)\de x \\
&\,\phantom{\frac{\de I_{\hom}(g,G)}{\de \cL^N}(x_0) \leq \limsup_{j\to\infty}\frac1{r_j^N}\liminf_{n\to\infty} \bigg(}+\int_{Q_{r_j}(x_0)\cap S_{u_{j,n}}} \!\!\!\!\!\psi\Big(\frac{x}{\eps_n},[u_{j,n}](x),\nu_{u_{j,n}}(x)\Big)\de \cH^{N-1}(x)\!\bigg) \\
\leq&\, \limsup_{j\to\infty}\liminf_{n\to\infty} \frac1{r_j^N}\bigg(\int_{Q_{r_j}(x_0)} \!\! W\Big(\frac{x}{\eps_n},\nabla u\Big(\frac{m_n k}{r_j}(x-x_0)\Big)+\nabla g(x_0)+G(x)-G(x_0) \! \Big)\de x \\
&  + \frac{r_j}{m_n k} \int_{Q_{r_j}(x_0)\cap \big(x_0+\frac{r_j}{m_n k}S_u\big)} \!\! \psi\Big(\frac{x}{\eps_n},[u]\Big(\frac{m_n k}{r_j}(x-x_0)\Big),\nu_{u}\Big(\frac{m_n k}{r_j}(x-x_0)\Big)\Big) \de \cH^{N-1}(x) \\
&  + \int_{Q_{r_j}(x_0)\cap S_g} \!\! \psi\Big(\frac{x}{\eps_n},[g](x),\nu_{g}(x)\Big) \de \cH^{N-1}(x)  + \int_{Q_{r_j}(x_0)\cap S_{h_j}} \!\! \psi\Big(\frac{x}{\eps_n},[h_j](x),\nu_{h_j}(x)\Big) \de \cH^{N-1}(x) \\
&  + \int_{Q_{r_j}(x_0)\cap S_{h_{j,n}}} \!\! \psi\Big(\frac{x}{\eps_n},[h_{j,n}](x),\nu_{h_{j,n}}(x)\Big) \de \cH^{N-1}(x)\bigg),
\end{split}
\end{equation*}
where we have used 
the subadditivity and the positive $1$-homogeneity of $\psi$ (see \ref{H4} and \ref{H5}) to obtain the second inequality.
Now, using, in order, \ref{H3}, the estimate in \eqref{817} and \eqref{818}, and finally \eqref{ubb09}, the last three integrals above vanish as first $n\to\infty$ and then $j\to\infty$.

We are left with one volume integral and one surface integral; by adding and subtracting $W(x/\eps_n,\nabla g(x_0)+\nabla u(m_nk(x-x_0)/r_j))$ in the volume integral and using \ref{H2}, H\"older's inequality, and \eqref{ubb09}, and
by changing variables according to 
\begin{equation}\label{ubb11}
Q_{r_j}(x_0)\ni x\mapsto \frac{m_n k}{r_j}(x-x_0)\eqqcolon z\in m_n kQ,
\end{equation} 
we have
\begin{equation*}
\begin{split}
&\,\frac{\de I_{\hom}(g,G)}{\de \cL^N}(x_0) \leq  \limsup_{j\to\infty}\liminf_{n\to\infty} \frac1{(m_n k)^N}\bigg(\int_{m_n kQ} \!\!W\Big(\frac{x_0}{\eps_n}+\frac{r_j z}{m_n k\eps_n},\nabla g(x_0)+\nabla u(z) \! \Big)\de z \\
&\,\qquad\qquad\qquad\qquad+\int_{{m_n kQ}\cap S_u} \psi\Big(\frac{x_0}{\eps_n}+\frac{r_j z}{m_n k\eps_n},[u](z),\nu_u(z)\Big)\de\cH^{N-1}(z)\bigg) \\
=&\, \limsup_{j\to\infty}\liminf_{n\to\infty} \frac1{(m_n k)^N}\bigg( \int_{m_nkQ} W\Big(\gamma_n+z+\frac1{m_n}\Big\langle\frac{r_j}{k\eps_n}\Big\rangle z,\nabla g(x_0)+\nabla u(z)\Big)\de z \\
&\,\qquad\qquad\qquad\qquad+\int_{m_n kQ\cap S_u} \psi\Big(\gamma_n+z+\frac1{m_n}\Big\langle\frac{r_j}{k\eps_n}\Big\rangle z,[u](z),\nu_u(z)\Big)\de\cH^{N-1}(z)\bigg) \\
=&\, \limsup_{j\to\infty}\liminf_{n\to\infty} \frac1{k^N}\bigg( \int_{kQ} W\Big(\gamma_n+z+\frac1{m_n}\Big\langle\frac{r_j}{k\eps_n}\Big\rangle z,\nabla g(x_0)+\nabla u(z)\Big)\de z \\
&\,\qquad\qquad\qquad\qquad+\int_{(kQ)\cap S_u} \psi\Big(\gamma_n+z+\frac1{m_n}\Big\langle\frac{r_j}{k\eps_n}\Big\rangle z,[u](z),\nu_u(z)\Big)\de\cH^{N-1}(z)\bigg),
\end{split}
\end{equation*}
where we have defined $\gamma_n\coloneqq \langle x_0/\eps_n\rangle$ and $m_n\coloneqq \lfloor r_j/k\eps_n\rfloor$, and used the decomposition
$$\frac1{m_n}\frac{r_j}{k\eps_n}z=\frac1{m_n}\Big(\Big\lfloor\frac{r_j}{k\eps_n}\Big\rfloor+\Big\langle\frac{r_j}{k\eps_n}\Big\rangle\Big)z=z+\frac1{m_n}\Big\langle\frac{r_j}{k\eps_n}\Big\rangle z$$
to get the first equality; the second equality follows from the $kQ$-periodicity of $u$ and from the $Q$-periodicity in the first variable of $W$ and $\psi$ (see \ref{H1}).

Upon noticing that $m_n^{-1}\langle r_j/k\eps_n\rangle\to0$ as $n\to\infty$, we can extract a subsequence $j\mapsto n(j)$ such that $\big|m_n^{-1}\langle r_j/k\eps_n\rangle z\big|<1/j$, so that, 
upon diagonalization and invoking \ref{H2c} and \ref{H3c}, we can write
\begin{equation*}
\begin{split}
&\,\frac{\de I_{\hom}(g,G)}{\de \cL^N}(x_0) \leq  \limsup_{j\to\infty} \frac1{k^N}\bigg( \int_{kQ} W(\gamma_{n(j)}+z,\nabla g(x_0)+\nabla u(z))\,\de z \\
&\,\qquad\qquad\qquad\qquad+\int_{(kQ)\cap S_u} \psi(\gamma_{n(j)}+z,[u](z),\nu_u(z))\de\cH^{N-1}(z)\bigg);
\end{split}
\end{equation*}
moreover, by using the definition of infimum in $H_{\hom}^{\tau}$ in \eqref{118b} (for $\tau=\gamma_{n(j)}$), both $k\in\N{}$ and $u\in SBV_{\#}(kQ;\R^d)$ can be chosen in such a way that 
\begin{equation*}
\begin{split}
&\,\frac{\de I_{\hom}(g,G)}{\de \cL^N}(x_0) \leq \limsup_{j\to\infty}\bigg( H_{\hom}^{\gamma_{n(j)}}(\nabla g(x_0),G(x_0)) +\frac1j\bigg)=H_{\hom}(\nabla g(x_0),G(x_0)),
\end{split}
\end{equation*}
where we have used the translation invariance property of $H_{\hom}$ (see Proposition~\ref{118a}) to obtain the last equality.
\qed

Putting \eqref{BEDLB} and \eqref{ubb00} together, we obtain that 
$$\frac{\de I_{\hom}(g,G)}{\de\cL^N}(x_0)=H_{\hom}(\nabla g(x_0),G(x_0))$$
for all the points $x_0\in\Omega$ satisfying the conditions stated at the beginning of Section~\ref{sect_bed}, thus proving the first part of the integral representation \eqref{106}.

\medskip

\subsection{The surface energy density}\label{sect_sed}
We tackle here the surface energy density $h_{\hom}$.
From now on, we consider a point $x_0\in S_g$. 
Recalling Proposition \ref{asprop2.22CF}, for every $U\in\cO(\Omega)$ and for every $(g,G)\in SD_p(U)$, the functional $U\mapsto I_{\hom}(g,G;U)$ in \eqref{105} is a measure. 
In particular, (see \eqref{growthcontrol}) there exists $C>0$ such that 
\begin{equation}\label{growthF}
I_{\hom}(g, G;U )\leq C\big(\cL^N(U)+ |D^s g|(U)\big). 
\end{equation} 
Observe that \eqref{growthF} guarantees that, for every $g \in 
SBV(\Omega;\R^d)$, the computation of the Radon--Nikod\'ym derivative $\displaystyle \frac{\de I_{\hom}(g, G)}{\de |D^s g|}(x_0)$ 
does not depend on $G$.
Indeed, let us consider $\{u_k\}\in\widetilde\cR_p(g,G;U)$ a recovery sequence for $I_{\hom}(g, G; U)$ and, by Theorems~\ref{Al} and~\ref{ctap}, let us consider $v\in SBV(U;\R^d)$ such that $\nabla v=-G$ and piece-wise constant functions $v_k\in SBV(U;\R^d)$ such that $v_k\to v$ in $L^1(U;\R^d)$. 
Finally, let us define $w_k\coloneqq u_k+v-v_k$, so that  $w_k \wpSD (g, 0)$ and therefore
$$I_{\hom}(g, 0;U)\leq \liminf_{k\to \infty} \bigg\{\int_U W\Big(\frac{x}{\varepsilon_k}, \nabla w_k(x)\Big)\de x+ \int_{U\cap S_{w_k}}\psi\Big(\frac{x}{\varepsilon_k}, [w_k](x),\nu_{w_k}(x)\Big)\de \cH^{N-1}(x)\bigg\}.$$
Thus, by invoking \ref{H2} and H\"{o}lder's inequality for the volume integrals, and first the sub-additivity of $\psi$ (see \ref{H5}) then the linear growth of $\psi$ (see \ref{H3}) for the surface integrals, we can estimate
\[\begin{split}
&I_{\hom}(g,0;U)-I_{\hom}(g, G;U)\leq \liminf_{k\to\infty} \bigg\{ \int_U \bigg(W\Big(\frac{x}{\varepsilon_k}, \nabla w_k(x)\Big)- W\Big(\frac{x}{\varepsilon_k}, \nabla u_k(x)\Big)\bigg)\de x \\
+& \int_{U \cap S_{w_k}}\psi\Big(\frac{x}{\varepsilon_k}, [w_k](x), \nu_{w_k}(x)\Big)\de\cH^{N-1}(x)-\int_{U \cap S_{u_k}}\psi\Big(\frac{x}{\varepsilon_k}, [u_k](x), \nu_{u_k}(x)\Big)\de \cH^{N-1}(x)\bigg\} \\
\leq&\,\liminf_{k\to \infty} C\bigg\{\int_U (1+ |G|^p(x))\,\de x + \int_{U \cap S_v}|[v](x)| \,\de\cH^{N-1}(x) +\int_{U \cap S_{v_k}}|[v_k](x)|\,\de \cH^{N-1}(x)\bigg\},
\end{split}\]
where $C>0$ is a suitable constant. By virtue of the estimate in \eqref{817} and by \eqref{818}, the two surface integrals in the last line above are bounded by the volume integral, so that, 
by exchanging the roles of $I_{\hom}(g,G; U)$ and $I_{\hom}(g, 0; U)$, we arrive at the conclusion that
$$|I_{\hom}(g,0;U)-I_{\hom}(g, G;U)|\leq C\int_U (1+ |G|^p(x))\,\de x,$$
for every $U \in \mathcal O(\Omega)$.
In turn, this guarantees that, for $\cH^{N-1}$-a.e.~$x_0 \in S_g$,
\begin{equation}\label{derFsur}
\frac{\de I_{\hom}(g,0)}{\de |D^s g|}(x_0)= \frac{\de I_{\hom}(g, G)}{\de |D^s g|}(x_0).
\end{equation}

In view of this, without loss of generality, we will consider $G=0$ for the rest of the proof.
The lower bound (see \eqref{SEDLB}) below will be obtained considering $g$ of the type $s_{\lambda,\nu}$ in \eqref{special}, with $(\lambda,\nu)\in(\R^d\setminus\{0\})\times\S^{n-1}$; the upper bound (see \eqref{SEDUB} below) will be obtained considering $g$ taking finitely many values, that is $g\in BV(\Omega;L)$ where $L\subset \R^d$ is a set with finite cardinality. In particular, the upper bound will also hold for functions of the type $g=s_{\lambda,\nu}$. 
To conclude, the general case will be obtained via standard approximation results as in \cite[Theorem 4.4, Step 2]{CF1997} (stemming from the ideas \cite[Proposition 4.8]{AMT}), so that this part of the proof (which relies on the continuity properties of $h_{\hom}$, see Proposition~\ref{gammahomcont}) will be omitted.

\subsubsection*{The surface energy density: lower bound}
In this section we prove that 
\begin{equation}\label{SEDLB}
\frac{\de I_{\hom}(g,0)}{\de |D^sg|}(x_0)\geq h_{\hom}([g](x_0),\nu_g(x_0)),
\end{equation}
by following the lines of \cite[Proposition 6.2]{BDV1996}.
Without loss of generality, we can suppose that $\nu_g(x_0)=e_1$ (the first vector of the canonical basis) and we denote $s_\lambda\coloneqq s_{\lambda,e_1}$, so that $\lambda=[g](x_0)$.
Let $\sigma \in (0,1)$ and define $Q_\sigma\coloneqq(-\sigma/2,\sigma/2) \times (-1/2,1/2)^{N-1}$. 
By the definition of relaxation in \eqref{105}, let $\{u_n\}\subset \cR_p(g,0;\Omega)$ 
be a recovery sequence such that $u_n \wpSD (s_\lambda,0)$ 
and
$$I_{\hom}(s_\lambda, 0; Q_\sigma)=\lim_{n\to \infty}\int_{Q_\sigma} W\Big(\frac{x}{\varepsilon_n}, \nabla u_n(x)\Big)\de x +\int_{Q_\sigma \cap S_{u_n}}\psi\Big(\frac{x}{\varepsilon_n}, [u_n](x), \nu_{u_n}(x)\Big)\de \cH^{N-1}(x).$$
We now substitute the sequence $\{u_n\}$ by a new sequence $\{\bar u_n\}\in SBV(\Omega;\R^d)\cap L^\infty(\Omega;\R^d)$ (this is possible thanks to Lemma~\ref{lemma2.20CF}) with the following properties: $\bar u_n \wpSD (s_\lambda,0)$; given $\{m_n\}$ a diverging sequence of integers such that $\beta_n\coloneqq m_n\eps_n\to0$ as $n\to\infty$, it holds
\begin{equation}\label{adessonumero}
\frac{1}{\beta_n^N} \int_{Q_\sigma} |\bar u_n(x) -s_\lambda(x)|\, \de x \to 0, \quad\text{and}\quad \frac{1}{\beta_n^N}\nabla \bar u_n \rightharpoonup 0\;\; \text{in $L^p(\Omega;\R^{d \times N})$}
\end{equation}
(the latter convergence is due to the metrizability of the weak convergence on bounded sets); and for every $\eta>0$,
\begin{equation}\label{Ishouldnumberit}
\begin{split}
I_{\hom}(s_\lambda, 0; Q_\sigma)+ \eta\geq\limsup_{n\to\infty} & \int_{Q_\sigma}  W\Big(\frac{x}{\varepsilon_n}, \nabla \bar u_n(x)\Big)\de x \\
&\,+ \int_{Q_\sigma \cap S_{\bar u_n}} \psi\Big(\frac{x}{\varepsilon_n}, [\bar u_n](x), \nu(\bar u_n)(x)\Big)\de \cH^{N-1}(x).
\end{split}
\end{equation}

For 
$\tau \in \{0\}\times \mathbb Z^{N-1}$, let $x_{n,\tau}\coloneqq\beta_n \tau$ and $Q_{n,\tau}\coloneqq x_{n,\tau}+ \beta_n Q_\sigma=\beta_n(Q_\sigma+\tau)$. Let $\tau(n)$ be the index corresponding to a 'minimal cube'  such that
\begin{equation}\label{cubi}
E_{\varepsilon_n}(\bar u_n; Q_{n, \tau(n)})\leq E_{\varepsilon_n} (\bar u_n; Q_{n,\tau}) 
\end{equation}
for every $\tau \in \{0\}\times \mathbb Z^{N-1}$ and $Q_{n,\tau} \subset Q_\sigma$. 
We now define $Q_{n}\coloneqq Q_{n,\tau(n)}$, $x_n\coloneqq x_{n,\tau(n)}$, and, for every $x \in Q_\sigma$, we let $w_n(x)\coloneqq \bar u_n (x_n+ \beta_n x)$. 
We claim that $w_n \in SBV(Q_\sigma;\R^d)\cap L^\infty(Q_\sigma;\R^d)$ and
\begin{equation}\label{w_n_prop}
\begin{cases}
\text{(i)} & \text{$\{w_n\}$ is equi-bounded;} \\
\text{(ii)} & \text{$w_n\to s_\lambda$  in $L^1(Q_\sigma; \R^d)$, as $n\to \infty$;} \\
\text{(iii)} & \text{$\displaystyle \int_{Q_\sigma} |\nabla w_n(x)|^p \, \de x \to 0$, as $n\to \infty$;} \\
\text{(iv)} & \!\! \begin{array}{l}
\text{$\displaystyle \limsup_{n\to \infty}\int_{S_{w_n}\cap Q_n} \psi\Big(\frac{x}{\alpha_n},[w_n](x), \nu_{w_n}(x)\Big)\de\cH^{N-1}(x) \leq I_{\hom}(s_\lambda,0; Q_\sigma)+ \eta$,}\\
\text{where $\{\alpha_n\}$ is a suitable vanishing sequence.}
\end{array}
\end{cases}
\end{equation} 
Indeed, \eqref{w_n_prop}(i) follows by construction and \eqref{w_n_prop}(ii) is obtained by changing variables according to $y=x_n+\beta_nx\in Q_n$, by the definition of $w_n$, by the choice of $\tau$, by the inclusion $Q_n\subset Q_\sigma$, and finally by the first limit in \eqref{adessonumero}.

In order to prove \eqref{w_n_prop}(iii), we observe that, by the boundedness of the energy, 
there exists a constant $C>0$ such that
\[\begin{split}
C \geq&\, E_{\varepsilon_n}(\bar u_n; Q_\sigma) \geq \sum_{\substack{\tau \in \{0\}\times \Z{N-1} \\ Q_{n,\tau} \subset Q_\sigma}} E_{\varepsilon_n}(\bar u_n; Q_{n,\tau}) \geq\Big\lfloor\frac{1}{\beta_n}\Big\rfloor^{N-1}E_{\varepsilon_n}(\bar u_n; Q_n) \\
\geq&\, \Big\lfloor\frac{1}{\beta_n}\Big\rfloor^{N-1}\int_{Q_n} W\Big(\frac{x}{\varepsilon_n},\nabla \bar u_n(x)\Big)\de x \geq \Big\lfloor\frac{1}{\beta_n}\Big\rfloor^{N-1}\bigg(C_W'\int_{Q_n} |\nabla \bar u_n(x)|^p\, \de x-c_W'\beta_n^N\sigma\bigg),
\end{split}\]
where the second inequality is due to the fact the cubes $Q_{n,\tau}$ are disjoint; the third inequality follows from counting them; in the fourth inequality we have used the non-negativity of $\psi$; in the last inequality we have exploited \ref{H8}.
Next, observe that (using the change of variables $y=x_n+\beta_nx\in Q_n$ and the inclusion $Q_n\subset Q_\sigma$ again)
\begin{align*}
\int_{Q_\sigma}|\nabla w_n(x)|^p\,\de x= \int_{Q_\sigma}|\beta_n \nabla \bar u_n (x_n+ \beta_n x)|^p\,\de x \leq \beta_n^{p-N} \int_{Q_\sigma}|\nabla \bar u_n(y)|^p\,\de y,
\end{align*}
(where in the second integrand we computed the gradient of the composed function),

whence
$$\int_{Q_\sigma}|\nabla w_n(x)|^p \de x \leq \frac{\beta_n^{p-N}}{C_W'}\bigg(C\Big\lfloor\frac1{\beta_n}\Big\rfloor^{1-N}+c_W'\beta_n^N\sigma\bigg)\to0\quad\text{as $n\to\infty$}.$$

We now prove \eqref{w_n_prop}(iv) with $\alpha_n= \varepsilon_n/\beta_n$. To this end, we observe that $x_n/\varepsilon_n= m_n\tau(n)\in \{0\}\times\Z{N-1}$ and so, by using the change of variables $y=x_n+\beta_nx\in Q_n$, the non-negativity of $W$, and the periodicity of $\psi$ (see \ref{H1}), we obtain
\[\begin{split}
&\int_{Q_\sigma \cap S_{w_n}} \psi\Big(\frac{x}{\alpha_n}, [w_n](x), \nu_{w_n}(x))\Big)\de \cH^{N-1}(x) \\
=& \int_{Q_\sigma\cap S_{w_n}} \psi\Big(\frac{\beta_n x}{\varepsilon_n}, [\bar u_n](x_n +\beta_n x), \nu_{w_n}(x_n+\beta_nx)\Big)\de \cH^{N-1}(x) \\
=&\, \frac{1}{\beta_n^{N-1}}\int_{x_n+ \beta_n (Q_\sigma\cap S_{w_n})}\psi\Big(\frac{y-x_n}{\varepsilon_n}, [\bar u_n](y), \nu_{\bar u_n}(y)\Big)\de \cH^{N-1}(y) \\
=&\, \frac{1}{\beta_n^{N-1}}\int_{Q_n\cap S_{\bar u_n}}\psi\Big(\frac{y}{\varepsilon_n},[\bar u_n](y), \nu_{\bar u_n}(y)\Big)\de \cH^{N-1}(y) \\
&\, \leq \frac{1}{\beta_n^{N-1}}E_{\varepsilon_n}(\bar u_n;Q_n)\leq \frac{1}{\beta_n^{N-1}}\Big\lfloor\frac{1}{\beta_n}\Big\rfloor^{1-N}E_{\varepsilon_n}(\bar u_n; Q_\sigma).
\end{split}\]
Thus \eqref{w_n_prop}(iv) follows from \eqref{Ishouldnumberit} since
$$\limsup_{n \to\infty}\int_{Q_\sigma\cap S_{u_n}} \psi\Big(\frac{x}{\alpha_n}, [w_n](x), \nu_{w_n}(x)\Big)\de \cH^{N-1}(x)\leq \limsup_{n\to\infty} E_{\varepsilon_n}(\bar u_n;Q_\sigma).$$

The sequence $\{w_n\}$ can now be modified into a new sequence $\{\tilde v_n\}$ such that $\nabla \tilde v_n=0$ a.e.~in $Q_\sigma$ as follows: for every $n\in\N$, we approximate via Theorem~\ref{app_thm} the pair $(0,-\nabla w_n)\in SD_p(Q_\sigma)$ by a sequence $\hat w_{n,k}$, so that $\hat v_{n,k}\coloneqq  w_n+\hat w_{n,k}\in SBV(Q_\sigma;\R^d)$ is such that
$$\lim_{k\to\infty} \hat v_{n,k}= w_n\quad\text{in $L^1(Q_\sigma;\R^d)$}\qquad\text{and}\qquad \nabla \hat v_{n,k}=0\quad\text{a.e.~in $Q_\sigma$;}$$
moreover, invoking \eqref{w_n_prop}(ii), we have that 
\begin{equation}\label{devocentrarla}
\lim_{n\to\infty}\lim_{k\to\infty} \hat v_{n,k}=s_\lambda\quad \text{in $L^1(Q_\sigma;\R^d)$.}
\end{equation}

Now,  by \ref{H3}, \ref{H5}, \eqref{817} and \eqref{818}, we have 
\[\begin{split}
&\int_{Q_\sigma\cap S_{\hat v_{n,k}}}  \psi\Big(\frac{x}{\alpha_n}, [\hat v_{n,k}](x), \nu_{\hat v_{n,k}}(x)\Big)\de \cH^{N-1}(x) \\
\leq&  \int_{Q_\sigma\cap S_{w_n}} \psi\Big(\frac{x}{\alpha_n}, [w_n](x), \nu_{w_n}(x)\Big)\de \cH^{N-1}(x) +C|D^s \hat w_{n,k}|(Q_\sigma)\\
\leq&  \int_{Q_\sigma\cap S_{w_n}} \psi\Big(\frac{x}{\alpha_n}, [w_n](x), \nu_{w_n}(x)\Big)\de \cH^{N-1}(x) +C\bigg( \|\nabla  w_n\|_{L^1(Q_\sigma;\R^{d\times N})} + \frac1k\bigg).
\end{split}\] 
Hence, by a standard diagonalization argument, by defining $\tilde v_n\coloneqq \hat v_{n, k(n)}$, we have that $\tilde v_n \wpSD (s_\lambda,0)$, 
 $\nabla \tilde v_n =0$, and
\begin{equation}\label{ladobbiamoinvocare}
\begin{split}
&\, \lim_{n\to \infty}\int_{Q_\sigma\cap S_{\tilde v_n}} \psi\Big(\frac{x}{\alpha_n}, [\tilde v_n](x), \nu_{\tilde v_n}(x)\Big)\de \cH^{N-1}(x) \\
\leq&\,  \liminf_{n \to \infty} \int_{Q_\sigma\cap S_{w_n}} \psi\Big(\frac{x}{\alpha_n}, [w_n](x), \nu_{w_n}(x)\Big)\de \cH^{N-1}(x).
\end{split}
\end{equation} 

The next step is to modify the sequence $\{\tilde v_n\}$ into a new sequence $\{v_n\}$ such that $v_n|_{\partial Q_\sigma}=s_\lambda|_{\partial Q_\sigma}$.
This can be achieved by defining the function
\begin{equation}\label{vndef}
	v_{n}\coloneqq \begin{cases}
s_\lambda  & \text{in $Q_\sigma \setminus (1-r_n)Q_\sigma$,} \\
\tilde v_n & \text{in $(1-r_n)Q_\sigma$,}
\end{cases}
\end{equation}
where $\{r_n\}\subset(0,1)$ is a sequence such that $\lim_{n \to\infty}r_n=1^-$ and, by \eqref{devocentrarla},
\begin{equation}\label{richia}
	\int_{\partial (1-r_n)Q_\sigma} |\tilde v_n(x)- s_\lambda(x) |\, \de \cH^{N-1}(x)<\frac{1}{n}.
	\end{equation}
Clearly $\nabla  v_n=0$ a.e.~in $Q_\sigma$ and, again by \eqref{devocentrarla}, $\lim_{n \to \infty} v_n =s_\lambda$ in $L^1(Q_\sigma;\R^d)$. Moreover, by \eqref{Lippsi}, \eqref{ladobbiamoinvocare}, and \eqref{w_n_prop}(iv), we have
\[\begin{split}
&\,\limsup_{n\to \infty}\int_{Q_\sigma\cap S_{v_n}}\psi\Big(\frac{x}{\alpha_n}, [v_n](x), \nu_{v_n}(x)\Big)\de \cH^{N-1}(x) \\
\leq &\, I_{\hom}(s_\lambda,0;Q_\sigma)+\eta+ \limsup_{n \to \infty}  |D^s \tilde v_n- D^s v_n|(\partial (1-r_n)Q_\sigma),
\end{split}\]
and the latter limit is $0$ by the choice of $r_n$, as consequence of \eqref{vndef} and \eqref{richia}. 

We conclude the proof by extending, without relabeling it, $v_n$ to the whole unit cube $Q$ by defining it as $s_\lambda$ in $Q\setminus Q_\sigma$, so that the previous inequality becomes 
$$\limsup_{n\to \infty}\int_{Q\cap S_{v_n}}\psi\Big(\frac{x}{\alpha_n}, [v_n](x), \nu_{v_n}(x)\Big)\de \cH^{N-1}(x)\leq  I_{\hom}(s_\lambda,0;Q)+\eta.$$
By Proposition~\ref{gammahomcont} and the change of variables $y=\alpha_n^{-1}x\in\alpha_n^{-1}Q$, the function $\bar v_n(y)\coloneqq v_n(\alpha_ny)$ belongs to $\cC^{\surface}(\lambda,e_1;\alpha_n^{-1}Q)$ (see \eqref{Csurface}), so that, recalling that we had set $\lambda=[g](x_0)$ and letting $\eta\to0$, we obtain \eqref{SEDLB}.
\qed

\subsubsection*{The surface energy density: upper bound}
In this section we prove that 
\begin{equation}\label{SEDUB}
\frac{\de I_{\hom}(g,0)}{\de |D^sg|}(x_0)\leq h_{\hom}([g](x_0),\nu_g(x_0))
\end{equation}
by following the lines of \cite[Proposition 6.2]{BDV1996}. Recall that by the preliminary discussion we made at the beginning of the section we will restrict ourselves to the case of piecewise constant functions $g$, that is $g\in BV(\Omega;L)$, where $L\subset\R^d$ has finite cardinality; naturally, such a function $g$ is also an element of $SBV(\Omega;\R^{d})$. We will obtain estimate \eqref{SEDUB} by using the abstract representation result contained in Theorem~\ref{thm31ABI}, for which we need to prove that our (localized) functional $I_{\hom}\colon SD_p(\Omega)\times\cO(\Omega)\to[0,+\infty)$ satisfies hypotheses \ref{thm31_a}--\ref{thm31_e} of Theorem~\ref{thm31ABI}.

As a consequence of \eqref{growthF}, for every $g \in BV(\Omega, L)$ and for every $U\in\cO(\Omega)$, the inequality $I_{\hom}(g, 0; U \cap S_g)\leq C|D^s g|(U \cap S_g)$ holds true, giving \ref{thm31_a}.
By Proposition~\ref{asprop2.22CF}, for every $g \in BV(\Omega; L)$, the set function $\cO(\Omega)\ni U\mapsto I_{\hom}(g, 0; U \cap S_g)$ is a measure, giving \ref{thm31_b}.
From the definition \eqref{105} of $I_{\hom}$ and from the locality property of the (sequence of) energies $\{E_{\eps_n}\}$, we obtain that $I_{\hom}(g, 0; U\cap S_g) = I_{\hom}(g_1,0; U\cap S_{g_1})$ whenever $g= g_1$ a.e.~in $U\in\cO(\Omega)$. Indeed, it suffices to notice that the competitors for $I_{\hom}(g,0)$ and $I_{\hom}(g_1,0)$ are the same. Therefore, condition \ref{thm31_c} is satisfied.
To show that condition \ref{thm31_d} holds, let us consider a sequence $\{g_n\}\subset SBV(\Omega;L)$ such that $g_n\to g$ pointwise a.e.; then, the fact that $L$ has finite cardinality entails that $g_n\to g$ in $L^1(\Omega;\R^d)$, and therefore that $g_n\wSD(g,0)\in SD(\Omega)$. Since $g\mapsto I_{\hom}(g,0;U)$ is lower semicontinuous for every $U\in\cO(\Omega)$ by definition of $\Gamma$-liminf, the desired inequality
$$I_{\hom}(g,0;U)\leq \liminf_{n \to \infty} I_{\hom}(g_n,0;U)$$
follows immediately. It remains to prove \ref{thm31_e}: as a matter of fact, we will prove a stronger condition, as it is obtained in the proof of \cite[Proposition~4.2]{BDV1996}. This translation invariance result, which is obtained following the argument in \cite[Lemma~3.7]{BDV1996}, provides then a sufficient condition for \ref{thm31_e}. 
We claim that, for every $z\in\R^N$ and every $U\in\cO(\Omega)$, we have
\begin{equation}\label{TI}
I_{\hom}(g, 0; U)= I_{\hom}(g(\cdot -z), 0;U+z).
\end{equation}
Indeed, let $z\in\R^N$ be given and observe that it can be approximated by means of a sequence of integers in the sense that there exists $\{z_n\}\subset\Z{N}$ such that $\eps_nz_n\to z$ as $n\to\infty$. Let now $U\in\cO(\Omega)$ be fixed, let $\{u_n\}\in\cR_p(g,0;U)$ be a recovery sequence for $I_{\hom}(g,0;U)$, and define $v_n\coloneqq u_n(\cdot-\varepsilon_n z_n)\colon U+z_n\to\R^d$.
Then we have 
\[\begin{split}
E_{\varepsilon_n}(u_n;U) = & \int_U W \Big(\frac{x + \varepsilon_n z_n}{\varepsilon_n}, \nabla u_n(x)\Big)\de x + \int_{U \cap S_{u_n}} \psi\Big(\frac{x + \varepsilon_n z_n}{\varepsilon_n}, [u_n](x), \nu_{u_n}(x)\Big) \de \cH^{N-1}(x) \\
=& \int_{U+\varepsilon_n z_n} W\Big(\frac{x}{\varepsilon_n}, \nabla v_n(x)\Big)\de x + \int_{(U+\varepsilon_n z_n)\cap S_{v_n}}
				\psi\Big( \frac{x}{\varepsilon_n}, [v_n](x), \nu_{v_n}(x)\Big)\de\cH^{N-1}(x).
\end{split}\]

Let now $V \subset \subset U$, so that, for $n$ sufficiently large we may assume $U + \varepsilon_n z_n \supseteq V + z$; hence, invoking the non-negativity of $W$ and $\psi$,
$$E_{\varepsilon_n}(u_n,U)\geq \int_{V+z}W\Big(\frac{x}{\varepsilon_n}, \nabla v_n(x)\Big)\de x+ \int_{(V+ z)\cap S_{v_n}} \psi\Big(\frac{x}{\varepsilon_n}, [v_n](x), \nu_{v_n}(x)\Big)\de\cH^{N-1}(x),$$
which yields  $I_{\hom}(g, 0;U) \geq I_{\hom}(g(\cdot -z), 0;V+z)$, since $v_n \wSD(g(\cdot-z), 0)$. 
By the arbitrariness of $V \subset \subset U$ we obtain that $I_{\hom}(g, 0; U)  \geq I_{\hom}(g(\cdot- z), 0;U+z)$. 
The reverse inequality can be obtained with the same reasoning, by defining $v_n\coloneqq u_n(\cdot +\varepsilon z_n)$. 
The translation invariance \eqref{TI} is proven, and this implies condition \ref{thm31_e}.

We are in position to apply Theorem~\ref{thm31ABI} and conclude that there exists a function $\tilde\psi\colon\Omega\times L\times L\times \S^{N-1}\to[0,+\infty)$ such that the integral representation 
$$I_{\hom}(g,0;U\cap S_g)=\int_{U \cap S_g} \tilde\psi \big(x,g^+(x),g^-(x),\nu_g(x)\big)\,\de \cH^{N-1}(x)$$ 
holds for every $g \in BV(\Omega;L)$ and for every $U\in \cO(\Omega)$.
Exactly with the same proof as in \cite[Lemma~3.7]{BDV1996} and \cite[Equation~(4.6)]{BDV1996}, one can prove that the density $\tilde\psi$ does not depend on the~$x$ variable and depends on $g$ only through its jump, so that there exists $\bar h_{\hom}\colon\R^d\times\S^{N-1}\to[0,+\infty)$ such that 
$$I_{\hom}(g,0;U\cap S_g)=\int_{U\cap S_g} \bar h_{\hom}([g](x),\nu_g(x))\,\de\cH^{N-1}(x).$$

On the other hand, we can get a precise estimate from above 
arguing as in \cite{BDV1996}.  
Upon defining the functional $J_{\hom}\colon BV(\Omega;L)\times\cO(\Omega)\to[0,+\infty)$ as
\begin{equation}\label{Jhom}
\!\! 
J_{\hom}(g;U)\coloneqq 
\inf\Big\{ \liminf_{n\to\infty} E_{\eps_n}(u_n;U): \{u_n\}\in\widetilde\cR_p(g,0;U),\, \sup_{n\in\N} \cH^{N-1}(U\cap S_{u_n}) < + \infty \Big\},
\end{equation}
we obtain
\begin{equation}\label{ineqsur}	
I_{\hom}(g,0;U\cap S_g)\leq J_{\hom}(g; U \cap S_g),
\end{equation}
for every  $g \in BV(\Omega; L)$ and every $U\in\cO(\Omega)$.
Now, \cite[Proposition 6.1]{BDV1996} grants that\footnote{We observe that in \cite{BDV1996} the standing hypothesis on the growth of the surface energy density $\psi$ is of the type $c(1+|\lambda|)\leq \psi(\lambda, \nu)\leq C(1+|\lambda|)$, in contrast with our \ref{H3}, the difference being that the measure $\cH^{N-1}(S_{u_n})$ cannot be controlled in our case. Nonetheless, this is circumvented by the introduction of the functional $J_{\hom}$ in \eqref{Jhom}, which provides an upper bound to $I_{\hom}$ and to which we can apply \cite[Proposition 6.1]{BDV1996}.}
$$J_{\hom}(g;U\cap S_g) \leq \int_{U\cap S_g} h_{\hom}([g](x),\nu_g(x))\,\de\cH^{N-1}(x),$$
where $h_{\hom}\colon\R^d\times \S^{N-1}\to[0,+\infty)$ is the functions defined in \eqref{108}; in turn, together with \eqref{ineqsur}, we obtain 
$$I_{\hom}(g, 0; U \cap S_g)\leq \int_{U \cap S_g} h_{\hom}([g](x), \nu_g(x))\,\de \cH^{N-1}(x),$$
whence 
$$\frac{\de I_{\hom}(g, 0)}{\de |D^s g|}(x_0) \leq h_{\hom}([g](x_0), \nu_g(x_0)) \qquad \text{for $\cH^{N-1}$-a.e.~$x_0 \in S_g$,}$$
which is \eqref{SEDUB} when $g
\in BV(\Omega, L)$.
Putting \eqref{SEDLB} and \eqref{SEDUB} together and keeping \eqref{derFsur} into account, we obtain that, for $g=s_{\lambda,\nu}\in BV(\Omega;L)$ and for all $G \in L^p(\Omega;\mathbb R^{d \times N})$, the equality
$$\frac{\de I_{\hom}(g,G)}{\de|D^sg|}(x_0)=h_{\hom}([g](x_0),\nu_g(x_0))$$
holds for all the points $x_0\in\Omega$ satisfying the conditions stated at the beginning of Section~\ref{sect_sed}.
To conclude, the equality in the general case, that is, for every $g\in SBV(\Omega;\R^d)$, is obtained via standard approximation results as in \cite[Theorem 4.4, Step 2]{CF1997} (stemming from the ideas \cite[Proposition 4.8]{AMT}), so that this part of the proof, which relies on the continuity properties of $h_{\hom}$ stated in Proposition~\ref{gammahomcont}, will be omitted.

Theorem~\ref{thm_main} is now completely proved. \qed


%
\bigskip 
\noindent\textbf{Acknowledgements.}
MA, MM, and EZ are members of the \emph{Gruppo Nazionale per l'Analisi Matematica, la Probabilit\`{a} e le loro Applicazioni} (GNAMPA) of the \emph{Instituto Nazionale di Alta Matematica ``F.~Severi''} (INdAM).
JM acknowledges support from FCT/Portugal through CAMGSD, IST-ID,
projects UIDB/04459/2020 and UIDP/04459/2020.
MM and EZ acknowledge funding from the GNAMPA Project 2020 \emph{Analisi variazionale di modelli non-locali nelle scienze applicate}.
MM is a member of the \emph{Integrated Additive Manufacturing} center at Politecnico di Torino.
MM acknowledges both that the present research fits in the scopes of the MIUR grant Dipartimenti di Eccellenza 2018-2022 (E11G18000350001) and partial support from the \emph{Starting grant per giovani ricercatori} of Politecnico di Torino.

\appendix
\section{Some technical results}\label{appa}
This appendix contains some technical results that are either reported here with no proof or proved for the reader's convenience since their proof is quite standard but in some measure different from the analogous results in the literature.

We start by showing that, for every $U \in \cO(\Omega)$ and for every $(g,G) \in SD_p(U)$, the localization $U\mapsto \tilde I_{\hom}(g, G; U)$ of the functional $\tilde I_{\hom}$ of \eqref{Ihomtilde}, defined as
\begin{equation}\label{calF}
\tilde I_{\hom}(g,G;U)\coloneqq \inf\Big\{\liminf_{n\to\infty} E_{\eps_n}(u_n): \{u_n\}\in \widetilde\cR_p(g,G;U)\Big\}
\end{equation}
(recalling the definition of the set $\widetilde\cR_p$ in \eqref{105-2}), is the trace of a Radon measure which is absolutely continuous with respect to $\cL^N + |D^s g|$.

To this end we recall a by now classical result of Fonseca and Mal\'y \cite{FMa}, which refines the De Giorgi--Letta criterion \cite{DGL1977} to establish sufficient conditions under which a functional is the restriction to open sets of a bounded Radon measure.
\begin{lemma}[{\cite{FMa}}]\label{FMaly}
Let $X$ be a locally compact Hausdorff space, let $\Pi \colon \cO(X) \to [0,+\infty]$ be a set function, and let $\mu$ be a finite Radon measure on $X$ satisfying
\begin{itemize}
\item[(i)] for every $U, V, Z \in \cO(X)$ such that $U \subset\subset V \subset \subset Z$, the following \emph{nested subadditivity} property holds: $\Pi(Z) \leq \Pi (V) + \Pi(Z\setminus \overline{U})$;
\item[(ii)] for every $U \in \cO(X)$ and for every $\eps > 0$ there exists $U_\eps \in  \cO(X)$ such that $U_\eps  \subset \subset U$ and $\Pi (U \setminus \overline{U_\eps})\leq \eps$;
\item[(iii)] $\Pi(X)\geq \mu(X)$;
\item[(iv)] for every $U \in \cO(X)$, it holds $\Pi(U) \leq \mu(\overline{U})$.
\end{itemize} 
Then $\Pi = \mu\res {\mathcal O}(X)$, that is, $\Pi$ is the restriction of the finite Radon measure $\mu$ to the open subsets of $X$.
\end{lemma}
The following proposition is in the spirit of \cite[Proposition~2.22]{CF1997}.
\begin{proposition}\label{asprop2.22CF}
Assume that \ref{H1}--\ref{H6} hold and let $(g, G) \in SD_p(\Omega)$. 
Then the localized functional $\cO(\Omega)\ni U\mapsto \tilde I_{\hom}(g, G;U)$ defined in \eqref{calF} is the trace on $\cO(\Omega)$  of a finite Radon measure on $\cB(\Omega)$. 
\end{proposition}
\begin{proof}  The proof relies on Lemma~ \ref{FMaly}: we will show that its hypotheses are satisfied by $\Pi(U)=\tilde I_{\hom}(g,G;U)$ and $\mu=\cL^N+|D^sg|$. In the rest of the proof, we will drop the tilde to keep the notation lighter.
First we prove that for every $U\in\cO(\Omega)$ and for every $(g, G)\in SD_p(\Omega)$ there exists a constant $C>0$ such that
\begin{equation}\label{growthcontrol}
I_{\hom}(g,G; U)\leq C\big(\cL^N(U)+ |D^s g|(U)\big). 
\end{equation}
We observe that by Theorem  \ref{app_thm} there exists $\{u_n\} \subset SBV(U;\mathbb R^d)$ such that $u_n\wpSD (g,G)$ and such that \eqref{appEST} and \eqref{centerline} hold.
%
%
Thus, the definition of $I_{\hom}(g, G; U)$, the linear growth condition \ref{H3} of $\psi$, and \eqref{pgrabove} entail that
\begin{equation}\label{boundmeas}
\begin{split}
&\,I_{\hom}(g, G ; U) \\
\leq&\, \liminf_{n\to \infty}\bigg\{\int_U W\Big(\frac{x}{\eps_n}, \nabla u_n(x)\Big)\de x + \int_{U \cap S_{u_n}} \psi\Big(\frac{x}{\eps_n}, [u_n](x), \nu_{u_n}(x)\Big)\de\cH^{N-1}(x)\bigg\} \\
\leq&\, \liminf_{n\to \infty}\bigg\{\int_U  C_W(1+|\nabla u_n(x)|^p) \de x + \int_{U \cap S_{u_n}} C_\psi |[u_n(x)]|\de \cH^{N-1}(x)\bigg\} \\
\leq&\, 
	C_W\Big(\cL^N(U)+\lVert G\rVert_{L^p(U;\R^{d\times N})}^p\Big)+C_\psi|Dg|(U)\eqqcolon\lambda(U),
\end{split}
\end{equation}
which implies \eqref{growthcontrol}.

We start proving condition (iv) in  Lemma~\ref{FMaly}. 
By the definition of $\Gamma$-limit , there exists a sequence $\{u_n\}\subset SBV(\Omega;\R^d)$ such that $u_n \wpSD(g,G)$ and along which 
$$I_{\hom}(g ,G; \Omega)=\lim_{n\to \infty }\bigg\{\int_{\Omega}W\Big(\frac{x}{\eps_n}, \nabla u_n(x)\Big)\de x+ \int_{\Omega \cap S_{u_n}}\psi\Big(\frac{x}{\eps_n}, [u_n](x), \nu_{u_n}(x)\Big)\de \cH^{N-1}(x)\bigg\}.$$
Upon the extraction of a subsequence, we know that 
$$W\Big(\frac{x}{\eps_n}, \nabla u_n(x)\Big) \de x + \psi\Big(\frac{x}{\eps_n}, [u_n](x), \nu_{u_n}(x)\Big)\de \cH^{N-1}\res S_{u_n} \wsto \mu\quad\text{in $\cM(\overline\Omega)$ as $n\to\infty$}$$ 
and 
\begin{equation}\label{2.22CF}
\mu({\overline \Omega})= I_{\hom}(g, G; \Omega) .
\end{equation}
On the other hand, for every $U \in \mathcal{O}(\Omega)$ we have that
\begin{equation}\label{2.23CF}
\begin{split}
I_{\hom}(g, G; U)\leq \liminf_{n\to \infty}\bigg\{&\int_U W\Big(\frac{x}{\eps_n},\nabla u_n(x)\Big) \de x \\
&+\int_{U \cap  S_{u_n} }\psi\Big(\frac{x}{\eps_n}, [u_n](x),\nu_{u_n}(x)\Big)\de \cH^{N-1}(x)\bigg\}\leq \mu({\overline U}).
\end{split}
\end{equation}

Next we prove condition (i) in Lemma \ref{FMaly}.
Consider $U,V,Z \in \cO(\Omega)$ such that $U\subset \subset V\subset \subset Z$. 
Fix $\eta >0$ and consider two sequences $\{u_n\}\subset SBV(V;\R^d)$ and $\{v_n\} \in SBV(Z\setminus \overline{U};\R^d)$ which are almost minimizing for $I_{\hom}$, that is,
\begin{equation*}
\begin{split}
&\lim_{n\to \infty}\bigg\{\int_V W\Big(\frac{x}{\varepsilon_n},\nabla u_n(x)\Big)\de x +\int_{V \cap S_ {u_n}}\psi\Big(\frac{x}{\varepsilon_n},[u_n](x), \nu_{u_n}(x)\Big) \de \cH^{N-1}(x)\bigg\} \\
\leq&\,  \eta + I_{\hom}(g, G; V), \\
&\lim_{n\to \infty}\bigg\{\int_{(Z\setminus \overline{U})} Z\Big(\frac{x}{\varepsilon_n},\nabla v_n(x)\Big)\de x +\int_{(Z\setminus \overline{U})\cap S_{ v_n}}\psi\Big(\frac{x}{\varepsilon_n},[v_n](x), \nu_{v_n}(x)\Big) \de\cH^{N-1}(x)\bigg\}\\
\leq&\, \eta + I_{\hom}(g, G; Z\setminus \overline{U}),
\end{split}
\end{equation*}
with $u_n \wpSD (g,G) \in SD_p(V)$, $v_n \wpSD  (g,G) \in SD_p(Z\setminus \overline{U})$.

In order to connect the functions without adding more interfaces, we argue as in the proof of \cite[Proposition 5.1]{MMZ}. 
For $\delta > 0$ small enough, consider
$$U_\delta\coloneqq \{x \in V: \dist{x}{\partial \overline{U}}< \delta\}.$$
For $x \in Z$, let $d(x) \coloneqq\dist{x}{U}$. Since the distance function to a fixed set is Lipschitz
continuous (see  \cite[Exercise 1.1]{Z}), we can apply the change of variables formula
(see \cite[Theorem 2, Section 3.4.3]{EG2015}), to obtain
$$\int_{U_\delta \setminus \overline{U}} |u_n(x)-v_n(x)| Jd(x)\,\de x=\int_0^\delta \bigg[\int_{d^{-1}(y)} |u_n(x)-v_n(x)|\,\de \cH^{N-1}(x)\bigg]\de y$$
and, since the Jacobian determinant $Jd(x)$ is bounded and $u_n - v_n \to 0$  in $L^1(V \cap (Z \setminus \overline{U});\R^d)$, it follows that for almost every $\varrho \in [0, \delta]$ we have
\begin{equation}\label{413}
\lim_{n\to \infty}\int_{d^{-1}(\varrho)}|u_n(x) - v_n(x)|\,\de\cH^{N-1}(x) = \lim_{n\to \infty}\int_{\partial U_\varrho} |u_n(x)-v_n(x)|\,\de \cH^{N-1}(x)=0.
\end{equation}
Fix $\varrho_0\in [0; \delta]$ such that \eqref{413} holds. We observe that $U_{\varrho_0}$ is a set with locally Lipschitz boundary since it is a level set of a Lipschitz function (see, \emph{e.g.}, \cite{EG2015}). Hence we can consider $u_n$ and $v_n $ on $\partial U_{\varrho_0}$ in the sense of traces and we can define
$$
w_n(x)\coloneqq\begin{cases}
u_n(x) & \text{if $x\in\overline{U}_{\varrho_0}$\,,}\\
v_n(x) & \text{if $x\in Z\setminus \overline{U}_{\varrho_0}$\,.}
\end{cases}
$$
By the choice of $\varrho_0$, the function $w_n$ is admissible for $I_{\hom}(g, G;Z)$; in particular 
$w_n \wpSD (g,G)  \in SD_p(Z;\R^d)$.
Thus we have
\begin{equation*}
\begin{split}
&I_{\hom}(g, G; Z)\\
\leq&\, \liminf_{n\to \infty}\bigg\{\int_Z W\Big(\frac{x}{\varepsilon_n}, \nabla w_n(x)\Big) \de x + \int_{W \cap S_ {w_n}}\psi\Big(\frac{x}{\varepsilon_n}, [w_n](x), \nu_{w_n}(x)\Big)\de\cH^{N-1}(x)\bigg\} \\
\leq&\, \liminf_{n\to \infty} \bigg\{\int_V W\Big (\frac{x}{\varepsilon_n}, \nabla u_n(x)\Big)\de x + \int_{V \cap S_{u_n}}\psi\Big(\frac{x}{\varepsilon_n},[u_n](x), \nu_{u_n}(x)\Big)\de\cH^{N-1}(x) \\
&\,\phantom{\liminf_{n\to \infty} \bigg\{} +\int_{Z\setminus \overline{U}} W\Big(\frac{x}{\varepsilon_n}, \nabla v_n(x)\Big)\de x + \int_{(Z\setminus\overline U) \cap S_{v_n}} \psi\Big(\frac{x}{\varepsilon_n}, [v_n](x), \nu_{v_n}(x)\Big)\de\cH^{N-1}(x) \bigg\}\\
&\,\phantom{\liminf_{n\to \infty} \bigg\{} + \limsup_{n\to \infty}\int_{ \partial U_{\varrho_0} \cap S_{w_n}}\psi\Big(\frac{x}{\varepsilon_n}, [w_n](x), \nu_{w_n}(x)\Big)\de \cH^{N-1}(x) \\
\leq&\, I_{\hom}(g, G; V)+ I_{\hom}(g, G; Z \setminus \overline{U})+2 \eta \\
&\,+ \limsup_{n\to \infty}\int_{ \partial U_{\varrho_0} \cap S_{w_n}}\psi\Big(\frac{x}{\varepsilon_n}, [w_n](x), \nu_{w_n}(x)\Big)\de\cH^{N-1}(x). 
	\end{split}
\end{equation*}
Observing that, by \ref{H3}, and \eqref{413}, the last surface integral converges to $0$, condition (i) follows by letting $\eta\to0^+$.

It remains to prove conditions (ii) and (iii) in Lemma \ref{FMaly}. 
To this end, fix $\eta >0$ and take $V\subset \subset Z$ such that $\mu (Z \setminus V)<\eta$. By (i), \eqref{2.22CF}, and \eqref{2.23CF}, it results
$$\mu(Z)\leq \mu(V)+\eta= \mu(\overline{\Omega})-\mu(\overline{\Omega}\setminus V)+\eta \leq I_{\hom}(g, G; \Omega)- I_{\hom}(g, G; \Omega \setminus \overline{V})+ \eta\leq I_{\hom}(g, G; Z)+\eta.$$
Letting $\eta \to 0^+$, we obtain $\mu(Z)\leq I_{\hom}(g, G; Z)$, which proves (ii).

Finally, 
fix $\eta>0$ and $Z\in\cO(\Omega)$ and take $K$ a compact set such that $K \subset \subset Z$ with $\lambda(Z\setminus K)<\eta$ (with $\lambda$ the measure in the right-hand side of \eqref{boundmeas}), and $V\in\cO(\Omega)$ an open set such that 
$K\subset \subset V\subset \subset Z$.
Using (i) and \eqref{2.23CF} we have
$$I_{\hom}(g, G; Z) \leq I_{\hom}(g, G; ZV)+ I_{\hom}(g, G; Z\setminus K) \leq \mu(\overline{V})+ \lambda(Z\setminus K)< \mu(Z)+\eta,$$
and (iii) follows by taking the limit $\eta \to 0^+$. 
The proposition is fully proved.
\end{proof}

We now report a result stating that the functional $I_{\hom}$ in \eqref{105} can be obtained along sequences that are bounded in $L^\infty$.
Let  $E_n$ be as in \eqref{100} and for every $g \in  L^\infty(\Omega;\mathbb R^d)\cap SBV(\Omega;\mathbb R^d)$, and $G \in L^p(\Omega;\mathbb R^{d \times N})$ define
\begin{equation}\label{Ihominfty}
I^\infty_{\hom}(g,G)\coloneqq \inf\Big\{\liminf_{n\to\infty} E_{\eps_n}(u_n): \{u_n\}\in \cR_p(g,G;\Omega) \hbox{ and }  \sup_{n\in\N{}}\lVert u_n\rVert_{L^\infty(\Omega;\R^d)}<+\infty \Big\}.
\end{equation}
Following the same arguments as in \cite[Lemma~2.20]{CF1997} one can prove that $I_{\hom}(g, G) = I^\infty_{\hom}(g,G)$ when $g \in L^\infty(\Omega;\mathbb R^d)\cap SBV(\Omega;\mathbb R^d)$, \emph{i.e.}, the additional $L^\infty$ bound on admissible sequences for \eqref{105} does not increase the energy, and it is used both in Sections~\ref{sect_bed} and~\ref{sect_sed}. 
\begin{lemma}\label{lemma2.20CF}
Let $p > 1$, $(g, G)\in \big(SBV(\Omega;\R^{d})\cap L^\infty(\Omega;\R^d)\big) \times L^p(\Omega;\R^{d \times N})$, and assume that Assumptions~\ref{ass} 
hold. Then
$$\tilde I_{\hom}(g,G)=I_{\hom}(g,G)= I^\infty_{\hom}(g,G).$$
\end{lemma}

Finally, 
 in the proof of the upper bound \eqref{SEDUB} for the surface energy density $h_{\hom}$, we used the integral representation result on partitions stated in \cite[Theorem 3.1]{ABI}, which we recall here for the reader's convenience.
\begin{theorem}\label{thm31ABI}
Let $L \subset \mathbb R^m$ be a subset with finite cardinality, and let $F\colon BV (\Omega;L) \times \cO(\Omega)\to [0, +\infty)$ be a functional satisfying the following conditions:
\begin{enumerate}
\item\label{thm31_a} there exists $\Lambda>0$ such that $0 \leq F(u; U) \leq  \Lambda \cH^{N-1}(U \cap S_u)$ for every $u \in BV(\Omega;L)$ and for every $U \in \cO(\Omega)$;
\item\label{thm31_b} $F(u;\cdot )$ is a measure for every $u \in BV (\Omega; L)$;
\item\label{thm31_c} $F(u; U) = F(v; U)$ whenever $u = v$ almost everywhere in $U\in\cO(\Omega)$;
\item\label{thm31_d} $u_h \to u$ a.e.~in $U$ implies that $F(u; U) \leq \liminf_{h \to +\infty}F(u_h;U)$ for every $U \in \cO(\Omega)$;
\item\label{thm31_e} for every $U \subset \subset \Omega$, there exists a continuous function $\omega_U\colon[0, +\infty) \to [0, +\infty)$ such that $\omega_U(0) = 0$ and $|F(u, V) - F(v, V+ z)|\leq \omega_U(|z|)\cH^{N-1}(U \cap S_u)$ whenever  $V \in \cO(U)$, $z \in \R^N$,  $|z|<\frac{\dist{U}{\partial \Omega}}{2}$, and $v(x + z) = u(x)$ in $V$.
\end{enumerate}
Then there exists a unique continuous function $f \colon\Omega \times L \times L\times  \S^{N-1}\to [0,\Lambda]$ such that $f(x, i, j, \nu) = f(x, j, i, -\nu)$, and the function $\displaystyle p \mapsto f\Big(x, i, j,\frac{p}{|p|}\Big) |p|$ is convex in $\R^N$ for every $x \in \Omega$, $i, j \in L$, and $F(u; U)$ is representable as
$$F(u; U) = \int_{U \cap S_u}f(x, u^+(x), u^-(x),  \nu_u(x))\, \de \cH^{N-1}(x),$$ 
for  every $u \in BV (\Omega; L)$ and for every $U \in \cO (\Omega)$.
\end{theorem}

\end{document}